\input amstex
%Sekretariats-Stildatei
%
%AMSPPT.STY modifiziert und erweitert von C. Krattenthaler

\def\next{AMS-SEKR}\ifx\styname\next \endinput\fi
\catcode`\@=11
\def\styname{AMS-SEKR}
\def\styversion{2.0}
{\W@{}\W@{\styname.STY - Version \styversion}\W@{}}
\hyphenation{acad-e-my acad-e-mies af-ter-thought anom-aly anom-alies
an-ti-deriv-a-tive an-tin-o-my an-tin-o-mies apoth-e-o-ses apoth-e-o-sis
ap-pen-dix ar-che-typ-al as-sign-a-ble as-sist-ant-ship as-ymp-tot-ic
asyn-chro-nous at-trib-uted at-trib-ut-able bank-rupt bank-rupt-cy
bi-dif-fer-en-tial blue-print busier busiest cat-a-stroph-ic
cat-a-stroph-i-cally con-gress cross-hatched data-base de-fin-i-tive
de-riv-a-tive dis-trib-ute dri-ver dri-vers eco-nom-ics econ-o-mist
elit-ist equi-vari-ant ex-quis-ite ex-tra-or-di-nary flow-chart
for-mi-da-ble forth-right friv-o-lous ge-o-des-ic ge-o-det-ic geo-met-ric
griev-ance griev-ous griev-ous-ly hexa-dec-i-mal ho-lo-no-my ho-mo-thetic
ideals idio-syn-crasy in-fin-ite-ly in-fin-i-tes-i-mal ir-rev-o-ca-ble
key-stroke lam-en-ta-ble light-weight mal-a-prop-ism man-u-script
mar-gin-al meta-bol-ic me-tab-o-lism meta-lan-guage me-trop-o-lis
met-ro-pol-i-tan mi-nut-est mol-e-cule mono-chrome mono-pole mo-nop-oly
mono-spline mo-not-o-nous mul-ti-fac-eted mul-ti-plic-able non-euclid-ean
non-iso-mor-phic non-smooth par-a-digm par-a-bol-ic pa-rab-o-loid
pa-ram-e-trize para-mount pen-ta-gon phe-nom-e-non post-script pre-am-ble
pro-ce-dur-al pro-hib-i-tive pro-hib-i-tive-ly pseu-do-dif-fer-en-tial
pseu-do-fi-nite pseu-do-nym qua-drat-ics quad-ra-ture qua-si-smooth
qua-si-sta-tion-ary qua-si-tri-an-gu-lar quin-tes-sence quin-tes-sen-tial
re-arrange-ment rec-tan-gle ret-ri-bu-tion retro-fit retro-fit-ted
right-eous right-eous-ness ro-bot ro-bot-ics sched-ul-ing se-mes-ter
semi-def-i-nite semi-ho-mo-thet-ic set-up se-vere-ly side-step sov-er-eign
spe-cious spher-oid spher-oid-al star-tling star-tling-ly
sta-tis-tics sto-chas-tic straight-est strange-ness strat-a-gem strong-hold
sum-ma-ble symp-to-matic syn-chro-nous topo-graph-i-cal tra-vers-a-ble
tra-ver-sal tra-ver-sals treach-ery turn-around un-at-tached un-err-ing-ly
white-space wide-spread wing-spread wretch-ed wretch-ed-ly Brown-ian
Eng-lish Euler-ian Feb-ru-ary Gauss-ian Grothen-dieck Hamil-ton-ian
Her-mit-ian Jan-u-ary Japan-ese Kor-te-weg Le-gendre Lip-schitz
Lip-schitz-ian Mar-kov-ian Noe-ther-ian No-vem-ber Rie-mann-ian
Schwarz-schild Sep-tem-ber
%Zustze (Sind auch in besp.exc einzutragen!):
form per-iods Uni-ver-si-ty cri-ti-sism for-ma-lism}
\Invalid@\nofrills
\Invalid@\usualspace
\newif\ifnofrills@
\def\nofrills@#1#2{\relaxnext@
  \DN@{\ifx\next\nofrills
    \nofrills@true\let#2\relax\DN@\nofrills{\nextii@}%
  \else
    \nofrills@false\def#2{#1}\let\next@\nextii@\fi
\next@}}
\def\usualspace@#1{\ifnofrills@\def\usualspace{#1}\fi}
\def\addto#1#2{\csname \expandafter\eat@\string#1@\endcsname
  \expandafter{\the\csname \expandafter\eat@\string#1@\endcsname#2}}
\newdimen\bigsize@
\def\big@#1#2{{\hbox{$\left#2\vcenter to#1\bigsize@{}%
  \right.\nulldelimiterspace\z@\m@th$}}}
\def\big{\big@\@ne}
\def\Big{\big@{1.5}}
\def\bigg{\big@\tw@}
\def\Bigg{\big@{2.5}}
\def\raggedcenter@{\leftskip\z@ plus.4\hsize \rightskip\leftskip
 \parfillskip\z@ \parindent\z@ \spaceskip.3333em \xspaceskip.5em
 \pretolerance9999\tolerance9999 \exhyphenpenalty\@M
 \hyphenpenalty\@M \let\\\linebreak}
\def\upperspecialchars{\def\ss{SS}\let\i=I\let\j=J\let\ae\AE\let\oe\OE
  \let\o\O\let\aa\AA\let\l\L}
\def\uppercasetext@#1{%
  {\spaceskip1.2\fontdimen2\the\font plus1.2\fontdimen3\the\font
   \upperspecialchars\uctext@#1$\m@th\aftergroup\eat@$}}
\def\uctext@#1$#2${\endash@#1-\endash@$#2$\uctext@}
\def\endash@#1-#2\endash@{\uppercase{#1}\if\notempty{#2}--\endash@#2\endash@\fi}
\def\runaway@#1{\DN@{#1}\ifx\envir@\next@
  \Err@{You seem to have a missing or misspelled \string\end#1 ...}%
  \let\envir@\empty\fi}
\newif\iftemp@
\def\notempty#1{TT\fi\def\test@{#1}\ifx\test@\empty\temp@false
  \else\temp@true\fi \iftemp@}
\font@\tensmc=cmcsc10
\font@\sevenex=cmex7
\font@\sevenit=cmti7
\font@\eightrm=cmr8 % preloaded in plain.tex
\font@\sixrm=cmr6 % preloaded in plain.tex
\font@\eighti=cmmi8     \skewchar\eighti='177 % preloaded
\font@\sixi=cmmi6       \skewchar\sixi='177   % preloaded
\font@\eightsy=cmsy8    \skewchar\eightsy='60 % preloaded
\font@\sixsy=cmsy6      \skewchar\sixsy='60   % preloaded
\font@\eightex=cmex8
\font@\eightbf=cmbx8 % preloaded in plain.tex
\font@\sixbf=cmbx6   % preloaded in plain.tex
\font@\eightit=cmti8 % preloaded in plain.tex
\font@\eightsl=cmsl8 % preloaded in plain.tex
\font@\eightsmc=cmcsc8
\font@\eighttt=cmtt8 % preloaded in plain.tex
%\font@\ninerm=cmr9
%\font@\ninei=cmmi9    \skewchar\ninei='177
%\font@\ninesy=cmsy9   \skewchar\ninesy='60
%\font@\nineex=cmex9
%\font@\ninebf=cmbx9
%\font@\nineit=cmti9
%\font@\ninesl=cmsl9
%\font@\ninesmc=cmcsc9
%\font@\ninemsa=msam9
%\font@\ninemsb=msbm9
%\font@\nineeufm=eufm9

%Ergnzung des fetten Small-Capitals-Fonts:
%\font@\eightbsmc=cmbcsc10 scaled 833
%\font@\tenbsmc=cmbcsc10
%\font@\twelvebsmc=cmbcsc10 scaled \magstep1
%\font@\fourteenbsmc=cmbcsc10 scaled \magstep2
%\font@\seventeenbsmc=cmbcsc10 scaled \magstep3
%\font@\twentybsmc=cmbcsc10 scaled \magstep4

\loadmsam
\loadmsbm
\loadeufm
\UseAMSsymbols
\newtoks\tenpoint@
\def\tenpoint{\normalbaselineskip12\p@
 \abovedisplayskip12\p@ plus3\p@ minus9\p@
 \belowdisplayskip\abovedisplayskip
 \abovedisplayshortskip\z@ plus3\p@
 \belowdisplayshortskip7\p@ plus3\p@ minus4\p@
 \textonlyfont@\rm\tenrm \textonlyfont@\it\tenit
 \textonlyfont@\sl\tensl \textonlyfont@\bf\tenbf
 \textonlyfont@\smc\tensmc \textonlyfont@\tt\tentt
%Ergnzung des fetten Small-Capitals-Fonts:
 \textonlyfont@\bsmc\tenbsmc
 \ifsyntax@ \def\big##1{{\hbox{$\left##1\right.$}}}%
  \let\Big\big \let\bigg\big \let\Bigg\big
 \else
  \textfont\z@=\tenrm  \scriptfont\z@=\sevenrm  \scriptscriptfont\z@=\fiverm
  \textfont\@ne=\teni  \scriptfont\@ne=\seveni  \scriptscriptfont\@ne=\fivei
  \textfont\tw@=\tensy \scriptfont\tw@=\sevensy \scriptscriptfont\tw@=\fivesy
  \textfont\thr@@=\tenex \scriptfont\thr@@=\sevenex
        \scriptscriptfont\thr@@=\sevenex
  \textfont\itfam=\tenit \scriptfont\itfam=\sevenit
        \scriptscriptfont\itfam=\sevenit
  \textfont\bffam=\tenbf \scriptfont\bffam=\sevenbf
        \scriptscriptfont\bffam=\fivebf
  \setbox\strutbox\hbox{\vrule height8.5\p@ depth3.5\p@ width\z@}%
  \setbox\strutbox@\hbox{\lower.5\normallineskiplimit\vbox{%
        \kern-\normallineskiplimit\copy\strutbox}}%
 \setbox\z@\vbox{\hbox{$($}\kern\z@}\bigsize@=1.2\ht\z@
 \fi
 \normalbaselines\rm\ex@.2326ex\jot3\ex@\the\tenpoint@}
\newtoks\eightpoint@
\def\eightpoint{\normalbaselineskip10\p@
 \abovedisplayskip10\p@ plus2.4\p@ minus7.2\p@
 \belowdisplayskip\abovedisplayskip
 \abovedisplayshortskip\z@ plus2.4\p@
 \belowdisplayshortskip5.6\p@ plus2.4\p@ minus3.2\p@
 \textonlyfont@\rm\eightrm \textonlyfont@\it\eightit
 \textonlyfont@\sl\eightsl \textonlyfont@\bf\eightbf
 \textonlyfont@\smc\eightsmc \textonlyfont@\tt\eighttt
%Ergnzung des fetten Small-Capitals-Fonts:
 \textonlyfont@\bsmc\eightbsmc
 \ifsyntax@\def\big##1{{\hbox{$\left##1\right.$}}}%
  \let\Big\big \let\bigg\big \let\Bigg\big
 \else
  \textfont\z@=\eightrm \scriptfont\z@=\sixrm \scriptscriptfont\z@=\fiverm
  \textfont\@ne=\eighti \scriptfont\@ne=\sixi \scriptscriptfont\@ne=\fivei
  \textfont\tw@=\eightsy \scriptfont\tw@=\sixsy \scriptscriptfont\tw@=\fivesy
  \textfont\thr@@=\eightex \scriptfont\thr@@=\sevenex
   \scriptscriptfont\thr@@=\sevenex
  \textfont\itfam=\eightit \scriptfont\itfam=\sevenit
   \scriptscriptfont\itfam=\sevenit
  \textfont\bffam=\eightbf \scriptfont\bffam=\sixbf
   \scriptscriptfont\bffam=\fivebf
 \setbox\strutbox\hbox{\vrule height7\p@ depth3\p@ width\z@}%
 \setbox\strutbox@\hbox{\raise.5\normallineskiplimit\vbox{%
   \kern-\normallineskiplimit\copy\strutbox}}%
 \setbox\z@\vbox{\hbox{$($}\kern\z@}\bigsize@=1.2\ht\z@
 \fi
 \normalbaselines\eightrm\ex@.2326ex\jot3\ex@\the\eightpoint@}

%Definition von 12-point, 14-point und 17-point Fonts
\font@\twelverm=cmr10 scaled\magstep1
\font@\twelveit=cmti10 scaled\magstep1
\font@\twelvesl=cmsl10 scaled\magstep1
\font@\twelvesmc=cmcsc10 scaled\magstep1
\font@\twelvett=cmtt10 scaled\magstep1
\font@\twelvebf=cmbx10 scaled\magstep1
\font@\twelvei=cmmi10 scaled\magstep1
\font@\twelvesy=cmsy10 scaled\magstep1
\font@\twelveex=cmex10 scaled\magstep1
\font@\twelvemsa=msam10 scaled\magstep1
\font@\twelveeufm=eufm10 scaled\magstep1
\font@\twelvemsb=msbm10 scaled\magstep1
\newtoks\twelvepoint@
\def\twelvepoint{\normalbaselineskip15\p@
 \abovedisplayskip15\p@ plus3.6\p@ minus10.8\p@
 \belowdisplayskip\abovedisplayskip
 \abovedisplayshortskip\z@ plus3.6\p@
 \belowdisplayshortskip8.4\p@ plus3.6\p@ minus4.8\p@
 \textonlyfont@\rm\twelverm \textonlyfont@\it\twelveit
 \textonlyfont@\sl\twelvesl \textonlyfont@\bf\twelvebf
 \textonlyfont@\smc\twelvesmc \textonlyfont@\tt\twelvett
%Ergnzung des fetten Small-Capitals-Fonts:
 \textonlyfont@\bsmc\twelvebsmc
 \ifsyntax@ \def\big##1{{\hbox{$\left##1\right.$}}}%
  \let\Big\big \let\bigg\big \let\Bigg\big
 \else
  \textfont\z@=\twelverm  \scriptfont\z@=\tenrm  \scriptscriptfont\z@=\sevenrm
  \textfont\@ne=\twelvei  \scriptfont\@ne=\teni  \scriptscriptfont\@ne=\seveni
  \textfont\tw@=\twelvesy \scriptfont\tw@=\tensy \scriptscriptfont\tw@=\sevensy
  \textfont\thr@@=\twelveex \scriptfont\thr@@=\tenex
        \scriptscriptfont\thr@@=\tenex
  \textfont\itfam=\twelveit \scriptfont\itfam=\tenit
        \scriptscriptfont\itfam=\tenit
  \textfont\bffam=\twelvebf \scriptfont\bffam=\tenbf
        \scriptscriptfont\bffam=\sevenbf
  \setbox\strutbox\hbox{\vrule height10.2\p@ depth4.2\p@ width\z@}%
  \setbox\strutbox@\hbox{\lower.6\normallineskiplimit\vbox{%
        \kern-\normallineskiplimit\copy\strutbox}}%
 \setbox\z@\vbox{\hbox{$($}\kern\z@}\bigsize@=1.4\ht\z@
 \fi
 \normalbaselines\rm\ex@.2326ex\jot3.6\ex@\the\twelvepoint@}

\def\headfonts{\twelvepoint\bf}

\font@\fourteenrm=cmr10 scaled\magstep2
\font@\fourteenit=cmti10 scaled\magstep2
\font@\fourteensl=cmsl10 scaled\magstep2
\font@\fourteensmc=cmcsc10 scaled\magstep2
\font@\fourteentt=cmtt10 scaled\magstep2
\font@\fourteenbf=cmbx10 scaled\magstep2
\font@\fourteeni=cmmi10 scaled\magstep2
\font@\fourteensy=cmsy10 scaled\magstep2
\font@\fourteenex=cmex10 scaled\magstep2
\font@\fourteenmsa=msam10 scaled\magstep2
\font@\fourteeneufm=eufm10 scaled\magstep2
\font@\fourteenmsb=msbm10 scaled\magstep2
\newtoks\fourteenpoint@
\def\fourteenpoint{\normalbaselineskip15\p@
 \abovedisplayskip18\p@ plus4.3\p@ minus12.9\p@
 \belowdisplayskip\abovedisplayskip
 \abovedisplayshortskip\z@ plus4.3\p@
 \belowdisplayshortskip10.1\p@ plus4.3\p@ minus5.8\p@
 \textonlyfont@\rm\fourteenrm \textonlyfont@\it\fourteenit
 \textonlyfont@\sl\fourteensl \textonlyfont@\bf\fourteenbf
 \textonlyfont@\smc\fourteensmc \textonlyfont@\tt\fourteentt
%Ergnzung des fetten Small-Capitals-Fonts:
 \textonlyfont@\bsmc\fourteenbsmc
 \ifsyntax@ \def\big##1{{\hbox{$\left##1\right.$}}}%
  \let\Big\big \let\bigg\big \let\Bigg\big
 \else
  \textfont\z@=\fourteenrm  \scriptfont\z@=\twelverm  \scriptscriptfont\z@=\tenrm
  \textfont\@ne=\fourteeni  \scriptfont\@ne=\twelvei  \scriptscriptfont\@ne=\teni
  \textfont\tw@=\fourteensy \scriptfont\tw@=\twelvesy \scriptscriptfont\tw@=\tensy
  \textfont\thr@@=\fourteenex \scriptfont\thr@@=\twelveex
        \scriptscriptfont\thr@@=\twelveex
  \textfont\itfam=\fourteenit \scriptfont\itfam=\twelveit
        \scriptscriptfont\itfam=\twelveit
  \textfont\bffam=\fourteenbf \scriptfont\bffam=\twelvebf
        \scriptscriptfont\bffam=\tenbf
  \setbox\strutbox\hbox{\vrule height12.2\p@ depth5\p@ width\z@}%
  \setbox\strutbox@\hbox{\lower.72\normallineskiplimit\vbox{%
        \kern-\normallineskiplimit\copy\strutbox}}%
 \setbox\z@\vbox{\hbox{$($}\kern\z@}\bigsize@=1.7\ht\z@
 \fi
 \normalbaselines\rm\ex@.2326ex\jot4.3\ex@\the\fourteenpoint@}

\def\chapheadfonts{\fourteenpoint\bf}

\font@\seventeenrm=cmr10 scaled\magstep3
\font@\seventeenit=cmti10 scaled\magstep3
\font@\seventeensl=cmsl10 scaled\magstep3
\font@\seventeensmc=cmcsc10 scaled\magstep3
\font@\seventeentt=cmtt10 scaled\magstep3
\font@\seventeenbf=cmbx10 scaled\magstep3
\font@\seventeeni=cmmi10 scaled\magstep3
\font@\seventeensy=cmsy10 scaled\magstep3
\font@\seventeenex=cmex10 scaled\magstep3
\font@\seventeenmsa=msam10 scaled\magstep3
\font@\seventeeneufm=eufm10 scaled\magstep3
\font@\seventeenmsb=msbm10 scaled\magstep3
\newtoks\seventeenpoint@
\def\seventeenpoint{\normalbaselineskip18\p@
 \abovedisplayskip21.6\p@ plus5.2\p@ minus15.4\p@
 \belowdisplayskip\abovedisplayskip
 \abovedisplayshortskip\z@ plus5.2\p@
 \belowdisplayshortskip12.1\p@ plus5.2\p@ minus7\p@
 \textonlyfont@\rm\seventeenrm \textonlyfont@\it\seventeenit
 \textonlyfont@\sl\seventeensl \textonlyfont@\bf\seventeenbf
 \textonlyfont@\smc\seventeensmc \textonlyfont@\tt\seventeentt
%Ergnzung des fetten Small-Capitals-Fonts:
 \textonlyfont@\bsmc\seventeenbsmc
 \ifsyntax@ \def\big##1{{\hbox{$\left##1\right.$}}}%
  \let\Big\big \let\bigg\big \let\Bigg\big
 \else
  \textfont\z@=\seventeenrm  \scriptfont\z@=\fourteenrm  \scriptscriptfont\z@=\twelverm
  \textfont\@ne=\seventeeni  \scriptfont\@ne=\fourteeni  \scriptscriptfont\@ne=\twelvei
  \textfont\tw@=\seventeensy \scriptfont\tw@=\fourteensy \scriptscriptfont\tw@=\twelvesy
  \textfont\thr@@=\seventeenex \scriptfont\thr@@=\fourteenex
        \scriptscriptfont\thr@@=\fourteenex
  \textfont\itfam=\seventeenit \scriptfont\itfam=\fourteenit
        \scriptscriptfont\itfam=\fourteenit
  \textfont\bffam=\seventeenbf \scriptfont\bffam=\fourteenbf
        \scriptscriptfont\bffam=\twelvebf
  \setbox\strutbox\hbox{\vrule height14.6\p@ depth6\p@ width\z@}%
  \setbox\strutbox@\hbox{\lower.86\normallineskiplimit\vbox{%
        \kern-\normallineskiplimit\copy\strutbox}}%
 \setbox\z@\vbox{\hbox{$($}\kern\z@}\bigsize@=2\ht\z@
 \fi
 \normalbaselines\rm\ex@.2326ex\jot5.2\ex@\the\seventeenpoint@}

\font@\rrrrrm=cmr10 scaled\magstep4
\font@\bigtitlefont=cmbx10 scaled\magstep4

\parindent1pc
\normallineskiplimit\p@
\newdimen\indenti \indenti=2pc
\def\pageheight#1{\vsize#1}
\def\pagewidth#1{\hsize#1%
   \captionwidth@\hsize \advance\captionwidth@-2\indenti}
\pagewidth{30pc} \pageheight{47pc}
\def\topmatter{%
 \ifx\undefined\msafam
 \else\font@\eightmsa=msam8 \font@\sixmsa=msam6
   \ifsyntax@\else \addto\tenpoint{\textfont\msafam=\tenmsa
              \scriptfont\msafam=\sevenmsa \scriptscriptfont\msafam=\fivemsa}%
     \addto\eightpoint{\textfont\msafam=\eightmsa \scriptfont\msafam=\sixmsa
              \scriptscriptfont\msafam=\fivemsa}%
   \fi
 \fi
 \ifx\undefined\msbfam
 \else\font@\eightmsb=msbm8 \font@\sixmsb=msbm6
   \ifsyntax@\else \addto\tenpoint{\textfont\msbfam=\tenmsb
         \scriptfont\msbfam=\sevenmsb \scriptscriptfont\msbfam=\fivemsb}%
     \addto\eightpoint{\textfont\msbfam=\eightmsb \scriptfont\msbfam=\sixmsb
         \scriptscriptfont\msbfam=\fivemsb}%
   \fi
 \fi
 \ifx\undefined\eufmfam
 \else \font@\eighteufm=eufm8 \font@\sixeufm=eufm6
   \ifsyntax@\else \addto\tenpoint{\textfont\eufmfam=\teneufm
       \scriptfont\eufmfam=\seveneufm \scriptscriptfont\eufmfam=\fiveeufm}%
     \addto\eightpoint{\textfont\eufmfam=\eighteufm
       \scriptfont\eufmfam=\sixeufm \scriptscriptfont\eufmfam=\fiveeufm}%
   \fi
 \fi
 \ifx\undefined\eufbfam
 \else \font@\eighteufb=eufb8 \font@\sixeufb=eufb6
   \ifsyntax@\else \addto\tenpoint{\textfont\eufbfam=\teneufb
      \scriptfont\eufbfam=\seveneufb \scriptscriptfont\eufbfam=\fiveeufb}%
    \addto\eightpoint{\textfont\eufbfam=\eighteufb
      \scriptfont\eufbfam=\sixeufb \scriptscriptfont\eufbfam=\fiveeufb}%
   \fi
 \fi
 \ifx\undefined\eusmfam
 \else \font@\eighteusm=eusm8 \font@\sixeusm=eusm6
   \ifsyntax@\else \addto\tenpoint{\textfont\eusmfam=\teneusm
       \scriptfont\eusmfam=\seveneusm \scriptscriptfont\eusmfam=\fiveeusm}%
     \addto\eightpoint{\textfont\eusmfam=\eighteusm
       \scriptfont\eusmfam=\sixeusm \scriptscriptfont\eusmfam=\fiveeusm}%
   \fi
 \fi
 \ifx\undefined\eusbfam
 \else \font@\eighteusb=eusb8 \font@\sixeusb=eusb6
   \ifsyntax@\else \addto\tenpoint{\textfont\eusbfam=\teneusb
       \scriptfont\eusbfam=\seveneusb \scriptscriptfont\eusbfam=\fiveeusb}%
     \addto\eightpoint{\textfont\eusbfam=\eighteusb
       \scriptfont\eusbfam=\sixeusb \scriptscriptfont\eusbfam=\fiveeusb}%
   \fi
 \fi
 \ifx\undefined\eurmfam
 \else \font@\eighteurm=eurm8 \font@\sixeurm=eurm6
   \ifsyntax@\else \addto\tenpoint{\textfont\eurmfam=\teneurm
       \scriptfont\eurmfam=\seveneurm \scriptscriptfont\eurmfam=\fiveeurm}%
     \addto\eightpoint{\textfont\eurmfam=\eighteurm
       \scriptfont\eurmfam=\sixeurm \scriptscriptfont\eurmfam=\fiveeurm}%
   \fi
 \fi
 \ifx\undefined\eurbfam
 \else \font@\eighteurb=eurb8 \font@\sixeurb=eurb6
   \ifsyntax@\else \addto\tenpoint{\textfont\eurbfam=\teneurb
       \scriptfont\eurbfam=\seveneurb \scriptscriptfont\eurbfam=\fiveeurb}%
    \addto\eightpoint{\textfont\eurbfam=\eighteurb
       \scriptfont\eurbfam=\sixeurb \scriptscriptfont\eurbfam=\fiveeurb}%
   \fi
 \fi
 \ifx\undefined\cmmibfam
 \else \font@\eightcmmib=cmmib8 \font@\sixcmmib=cmmib6
   \ifsyntax@\else \addto\tenpoint{\textfont\cmmibfam=\tencmmib
       \scriptfont\cmmibfam=\sevencmmib \scriptscriptfont\cmmibfam=\fivecmmib}%
    \addto\eightpoint{\textfont\cmmibfam=\eightcmmib
       \scriptfont\cmmibfam=\sixcmmib \scriptscriptfont\cmmibfam=\fivecmmib}%
   \fi
 \fi
 \ifx\undefined\cmbsyfam
 \else \font@\eightcmbsy=cmbsy8 \font@\sixcmbsy=cmbsy6
   \ifsyntax@\else \addto\tenpoint{\textfont\cmbsyfam=\tencmbsy
      \scriptfont\cmbsyfam=\sevencmbsy \scriptscriptfont\cmbsyfam=\fivecmbsy}%
    \addto\eightpoint{\textfont\cmbsyfam=\eightcmbsy
      \scriptfont\cmbsyfam=\sixcmbsy \scriptscriptfont\cmbsyfam=\fivecmbsy}%
   \fi
 \fi
 \let\topmatter\relax}
\def\chapterno@{\uppercase\expandafter{\romannumeral\chaptercount@}}
\newcount\chaptercount@
\def\chapter{\nofrills@{\afterassignment\chapterno@
                        CHAPTER \global\chaptercount@=}\chapter@
 \DNii@##1{\leavevmode\hskip-\leftskip
   \rlap{\vbox to\z@{\vss\centerline{\eightpoint
   \chapter@##1\unskip}\baselineskip2pc\null}}\hskip\leftskip
   \nofrills@false}%
 \FN@\next@}
\newbox\titlebox@

%Uppercase ist abgestellt.
\def\title{\nofrills@{\relax}\title@%
 \DNii@##1\endtitle{\global\setbox\titlebox@\vtop{\tenpoint\bf
 \raggedcenter@\ignorespaces
 \baselineskip1.3\baselineskip\title@{##1}\endgraf}%
 \ifmonograph@ \edef\next{\the\leftheadtoks}\ifx\next\empty
    \leftheadtext{##1}\fi
 \fi
 \edef\next{\the\rightheadtoks}\ifx\next\empty \rightheadtext{##1}\fi
 }\FN@\next@}
\newbox\authorbox@
\def\author#1\endauthor{\global\setbox\authorbox@
 \vbox{\tenpoint\smc\raggedcenter@\ignorespaces
 #1\endgraf}\relaxnext@ \edef\next{\the\leftheadtoks}%
 \ifx\next\empty\leftheadtext{#1}\fi}
\newbox\affilbox@
\def\affil#1\endaffil{\global\setbox\affilbox@
 \vbox{\tenpoint\raggedcenter@\ignorespaces#1\endgraf}}
\newcount\addresscount@
\addresscount@\z@
\def\address#1\endaddress{\global\advance\addresscount@\@ne
  \expandafter\gdef\csname address\number\addresscount@\endcsname
  {\vskip12\p@ minus6\p@\noindent\eightpoint\smc\ignorespaces#1\par}}
\def\email{\nofrills@{\eightpoint{\it E-mail\/}:\enspace}\email@
  \DNii@##1\endemail{%
  \expandafter\gdef\csname email\number\addresscount@\endcsname
  {\def\usualspace{{\it\enspace}}\smallskip\noindent\eightpoint\email@
  \ignorespaces##1\par}}%
 \FN@\next@}
\def\thedate@{}
\def\date#1\enddate{\gdef\thedate@{\tenpoint\ignorespaces#1\unskip}}
\def\thethanks@{}
\def\thanks#1\endthanks{\gdef\thethanks@{\eightpoint\ignorespaces#1.\unskip}}
\def\thekeywords@{}
\def\keywords{\nofrills@{{\it Key words and phrases.\enspace}}\keywords@
 \DNii@##1\endkeywords{\def\thekeywords@{\def\usualspace{{\it\enspace}}%
 \eightpoint\keywords@\ignorespaces##1\unskip.}}%
 \FN@\next@}
\def\thesubjclass@{}
\def\subjclass{\nofrills@{{\rm2020 {\it Mathematics Subject
   Classification\/}.\enspace}}\subjclass@
 \DNii@##1\endsubjclass{\def\thesubjclass@{\def\usualspace
  {{\rm\enspace}}\eightpoint\subjclass@\ignorespaces##1\unskip.}}%
 \FN@\next@}
\newbox\abstractbox@
\def\abstract{\nofrills@{{\smc Abstract.\enspace}}\abstract@
 \DNii@{\setbox\abstractbox@\vbox\bgroup\noindent$$\vbox\bgroup
  \def\envir@{abstract}\advance\hsize-2\indenti
  \usualspace@{{\enspace}}\eightpoint \noindent\abstract@\ignorespaces}%
 \FN@\next@}
\def\endabstract{\par\unskip\egroup$$\egroup}
\def\widestnumber#1#2{\begingroup\let\head\null\let\subhead\empty
   \let\subsubhead\subhead
   \ifx#1\head\global\setbox\tocheadbox@\hbox{#2.\enspace}%
   \else\ifx#1\subhead\global\setbox\tocsubheadbox@\hbox{#2.\enspace}%
   \else\ifx#1\key\bgroup\let\endrefitem@\egroup
        \key#2\endrefitem@\global\refindentwd\wd\keybox@
   \else\ifx#1\no\bgroup\let\endrefitem@\egroup
        \no#2\endrefitem@\global\refindentwd\wd\nobox@
   \else\ifx#1\page\global\setbox\pagesbox@\hbox{\quad\bf#2}%
   \else\ifx#1\item\setboxz@h{#2}\global\rosteritemwd\wdz@
        \global\advance\rosteritemwd by.5\parindent
   \else\message{\string\widestnumber is not defined for this option
   (\string#1)}%
\fi\fi\fi\fi\fi\fi\endgroup}
\newif\ifmonograph@
\def\Monograph{\monograph@true \let\headmark\rightheadtext
  \let\varindent@\indent \def\headfont@{\bf}\def\proclaimheadfont@{\smc}%
  \def\demofont@{\smc}}
%Bei Proclaim,...: Einrcken.
\let\varindent@\indent

\newbox\tocheadbox@    \newbox\tocsubheadbox@
\newbox\tocbox@
\def\toc{\toc@{Contents}}
\def\newtocdefs{%
   \def \title##1\endtitle
       {\penaltyandskip@\z@\smallskipamount
        \hangindent\wd\tocheadbox@\noindent{\bf##1}}%
   \def \chapter##1{%
        Chapter \uppercase\expandafter{\romannumeral##1.\unskip}\enspace}%
   \def \specialhead##1\endspecialhead
       {\par\hangindent\wd\tocheadbox@ \noindent##1\par}%
   \def \head##1 ##2\endhead
       {\par\hangindent\wd\tocheadbox@ \noindent
        \if\notempty{##1}\hbox to\wd\tocheadbox@{\hfil##1\enspace}\fi
        ##2\par}%
   \def \subhead##1 ##2\endsubhead
       {\par\vskip-\parskip {\normalbaselines
        \advance\leftskip\wd\tocheadbox@
        \hangindent\wd\tocsubheadbox@ \noindent
        \if\notempty{##1}\hbox to\wd\tocsubheadbox@{##1\unskip\hfil}\fi
         ##2\par}}%
   \def \subsubhead##1 ##2\endsubsubhead
       {\par\vskip-\parskip {\normalbaselines
        \advance\leftskip\wd\tocheadbox@
        \hangindent\wd\tocsubheadbox@ \noindent
        \if\notempty{##1}\hbox to\wd\tocsubheadbox@{##1\unskip\hfil}\fi
        ##2\par}}}
\def\toc@#1{\relaxnext@
   \def\page##1%
       {\unskip\penalty0\null\hfil
        \rlap{\hbox to\wd\pagesbox@{\quad\hfil##1}}\hfilneg\penalty\@M}%
 \DN@{\ifx\next\nofrills\DN@\nofrills{\nextii@}%
      \else\DN@{\nextii@{{#1}}}\fi
      \next@}%
 \DNii@##1{%
\ifmonograph@\bgroup\else\setbox\tocbox@\vbox\bgroup
   \centerline{\headfont@\ignorespaces##1\unskip}\nobreak
   \vskip\belowheadskip \fi
   \setbox\tocheadbox@\hbox{0.\enspace}%
   \setbox\tocsubheadbox@\hbox{0.0.\enspace}%
   \leftskip\indenti \rightskip\leftskip
   \setbox\pagesbox@\hbox{\bf\quad000}\advance\rightskip\wd\pagesbox@
   \newtocdefs
 }%
 \FN@\next@}
\def\endtoc{\par\egroup}
\let\pretitle\relax
\let\preauthor\relax
\let\preaffil\relax
\let\predate\relax
\let\preabstract\relax
\let\prepaper\relax
\def\dedicatory #1\enddedicatory{\def\preabstract{{\medskip
  \eightpoint\it \raggedcenter@#1\endgraf}}}
\def\thetranslator@{}
\def\translator#1\endtranslator{\def\thetranslator@{\nobreak\medskip
 \line{\eightpoint\hfil Translated by \uppercase{#1}\qquad\qquad}\nobreak}}
\outer\def\endtopmatter{\runaway@{abstract}%
 \edef\next{\the\leftheadtoks}\ifx\next\empty
  \expandafter\leftheadtext\expandafter{\the\rightheadtoks}\fi
 \ifmonograph@\else
   \ifx\thesubjclass@\empty\else \makefootnote@{}{\thesubjclass@}\fi
   \ifx\thekeywords@\empty\else \makefootnote@{}{\thekeywords@}\fi
   \ifx\thethanks@\empty\else \makefootnote@{}{\thethanks@}\fi
 \fi
  \pretitle
  \ifmonograph@ \topskip7pc \else \topskip4pc \fi
  \box\titlebox@
  \topskip10pt% reset to normal value
  \preauthor
  \ifvoid\authorbox@\else \vskip2.5pc plus1pc \unvbox\authorbox@\fi
  \preaffil
  \ifvoid\affilbox@\else \vskip1pc plus.5pc \unvbox\affilbox@\fi
  \predate
  \ifx\thedate@\empty\else \vskip1pc plus.5pc \line{\hfil\thedate@\hfil}\fi
  \preabstract
  \ifvoid\abstractbox@\else \vskip1.5pc plus.5pc \unvbox\abstractbox@ \fi
  \ifvoid\tocbox@\else\vskip1.5pc plus.5pc \unvbox\tocbox@\fi
  \prepaper
  \vskip2pc plus1pc
}
\def\document{\let\fontlist@\relax\let\alloclist@\relax
  \tenpoint}

%Modifizierte Head-Skips
\newskip\aboveheadskip       \aboveheadskip1.8\bigskipamount
\newdimen\belowheadskip      \belowheadskip1.8\medskipamount

\def\headfont@{\smc}
\def\penaltyandskip@#1#2{\relax\ifdim\lastskip<#2\relax\removelastskip
      \ifnum#1=\z@\else\penalty@#1\relax\fi\vskip#2%
  \else\ifnum#1=\z@\else\penalty@#1\relax\fi\fi}
\def\nobreak{\penalty\@M
  \ifvmode\def\penalty@{\let\penalty@\penalty\count@@@}%
  \everypar{\let\penalty@\penalty\everypar{}}\fi}
\let\penalty@\penalty
\def\heading#1\endheading{\head#1\endhead}

\def\specialheadfont@{\bf}
\outer\def\specialhead{\par\penaltyandskip@{-200}\aboveheadskip
  \begingroup\interlinepenalty\@M\rightskip\z@ plus\hsize \let\\\linebreak
  \specialheadfont@\noindent\ignorespaces}
\def\endspecialhead{\par\endgroup\nobreak\vskip\belowheadskip}
%\outer\def\head#1\endhead{\par\penaltyandskip@{-200}\aboveheadskip
% {\headfont@\raggedcenter@\interlinepenalty\@M
% \ignorespaces#1\endgraf}\nobreak
% \vskip\belowheadskip
% \headmark{#1}}
\let\headmark\eat@
\newskip\subheadskip       \subheadskip\medskipamount
\def\subheadfont@{\bf}
\outer\def\subhead{\nofrills@{.\enspace}\subhead@
 \DNii@##1\endsubhead{\par\penaltyandskip@{-100}\subheadskip
  \varindent@{\usualspace@{{\subheadfont@\enspace}}%
 \subheadfont@\ignorespaces##1\unskip\subhead@}\ignorespaces}%
 \FN@\next@}
\outer\def\subsubhead{\nofrills@{.\enspace}\subsubhead@
 \DNii@##1\endsubsubhead{\par\penaltyandskip@{-50}\medskipamount
      {\usualspace@{{\it\enspace}}%
  \it\ignorespaces##1\unskip\subsubhead@}\ignorespaces}%
 \FN@\next@}
\def\proclaimheadfont@{\bf}
\outer\def\proclaim{\runaway@{proclaim}\def\envir@{proclaim}%
  \nofrills@{.\enspace}\proclaim@
 \DNii@##1{\penaltyandskip@{-100}\medskipamount\varindent@
   \usualspace@{{\proclaimheadfont@\enspace}}\proclaimheadfont@
   \ignorespaces##1\unskip\proclaim@
  \sl\ignorespaces}% 
 \FN@\next@}
\outer\def\endproclaim{\let\envir@\relax\par\rm
  \penaltyandskip@{55}\medskipamount}
\def\demoheadfont@{\it}
\def\demo{\runaway@{proclaim}\nofrills@{.\enspace}\demo@
     \DNii@##1{\par\penaltyandskip@\z@\medskipamount
  {\usualspace@{{\demoheadfont@\enspace}}%
  \varindent@\demoheadfont@\ignorespaces##1\unskip\demo@}\rm
  \ignorespaces}\FN@\next@}
\def\enddemo{\par\medskip}
\def\qed{\ifhmode\unskip\nobreak\fi\quad\ifmmode\square\else$\m@th\square$\fi}
\let\remark\demo
%\endremark=\enddemo
\let\endremark\enddemo
%Headfont fr Definition wie bei Beweisen.
\def\definition{\runaway@{proclaim}%
  \nofrills@{.\demoheadfont@\enspace}\definition@
        \DNii@##1{\penaltyandskip@{-100}\medskipamount
        {\usualspace@{{\demoheadfont@\enspace}}%
        \varindent@\demoheadfont@\ignorespaces##1\unskip\definition@}%
        \rm \ignorespaces}\FN@\next@}

%\let\example\definition
%\let\endexample\enddefinition
%Modifikation:

\newdimen\rosteritemwd
\newcount\rostercount@
\newif\iffirstitem@
\let\plainitem@\item
\newtoks\everypartoks@
\def\par@{\everypartoks@\expandafter{\the\everypar}\everypar{}}
\def\roster{\edef\leftskip@{\leftskip\the\leftskip}%
 \relaxnext@
 \rostercount@\z@  
 \def\item{\FN@\rosteritem@}% 
 \DN@{\ifx\next\runinitem\let\next@\nextii@\else
  \let\next@\nextiii@\fi\next@}%
 \DNii@\runinitem  
  {\unskip  
   \DN@{\ifx\next[\let\next@\nextii@\else
    \ifx\next"\let\next@\nextiii@\else\let\next@\nextiv@\fi\fi\next@}%
   \DNii@[####1]{\rostercount@####1\relax
    \enspace{\rm(\number\rostercount@)}~\ignorespaces}%
   \def\nextiii@"####1"{\enspace{\rm####1}~\ignorespaces}%
   \def\nextiv@{\enspace{\rm(1)}\rostercount@\@ne~}%
   \par@\firstitem@false  
   \FN@\next@}% 
 \def\nextiii@{\par\par@  
  \penalty\@m\smallskip\vskip-\parskip
  \firstitem@true}%  
 \FN@\next@}
\def\rosteritem@{\iffirstitem@\firstitem@false\else\par\vskip-\parskip\fi
 \leftskip3\parindent\noindent  
 \DNii@[##1]{\rostercount@##1\relax
  \llap{\hbox to2.5\parindent{\hss\rm(\number\rostercount@)}%
   \hskip.5\parindent}\ignorespaces}%
 \def\nextiii@"##1"{%
  \llap{\hbox to2.5\parindent{\hss\rm##1}\hskip.5\parindent}\ignorespaces}%
 \def\nextiv@{\advance\rostercount@\@ne
  \llap{\hbox to2.5\parindent{\hss\rm(\number\rostercount@)}%
   \hskip.5\parindent}}%
 \ifx\next[\let\next@\nextii@\else\ifx\next"\let\next@\nextiii@\else
  \let\next@\nextiv@\fi\fi\next@}

\newif\ifnextRunin@
\def\endroster{\relaxnext@
 \par\leftskip@  
 \penalty-50 \vskip-\parskip\smallskip  
 \DN@{\ifx\next\Runinitem\let\next@\relax
  \else\nextRunin@false\let\item\plainitem@  
   \ifx\next\par 
    \DN@\par{\everypar\expandafter{\the\everypartoks@}}%
   \else  
    \DN@{\noindent\everypar\expandafter{\the\everypartoks@}}%
  \fi\fi\next@}%
 \FN@\next@}
\newcount\rosterhangafter@
\def\Runinitem#1\roster\runinitem{\relaxnext@
 \rostercount@\z@ 
 \def\item{\FN@\rosteritem@}%  
 \def\runinitem@{#1}%  
 \DN@{\ifx\next[\let\next\nextii@\else\ifx\next"\let\next\nextiii@
  \else\let\next\nextiv@\fi\fi\next}%
 \DNii@[##1]{\rostercount@##1\relax
  \def\item@{{\rm(\number\rostercount@)}}\nextv@}%
 \def\nextiii@"##1"{\def\item@{{\rm##1}}\nextv@}%
 \def\nextiv@{\advance\rostercount@\@ne
  \def\item@{{\rm(\number\rostercount@)}}\nextv@}%
 \def\nextv@{\setbox\z@\vbox  
  {\ifnextRunin@\noindent\fi  
  \runinitem@\unskip\enspace\item@~\par  
  \global\rosterhangafter@\prevgraf}% 
  \firstitem@false  
  \ifnextRunin@\else\par\fi  
  \hangafter\rosterhangafter@\hangindent3\parindent
  \ifnextRunin@\noindent\fi  
  \runinitem@\unskip\enspace 
  \item@~\ifnextRunin@\else\par@\fi  
  \nextRunin@true\ignorespaces}%  
 \FN@\next@}
\def\footmarkform@#1{$\m@th^{#1}$}
\let\thefootnotemark\footmarkform@
\def\makefootnote@#1#2{\insert\footins
 {\interlinepenalty\interfootnotelinepenalty
 \eightpoint\splittopskip\ht\strutbox\splitmaxdepth\dp\strutbox
 \floatingpenalty\@MM\leftskip\z@\rightskip\z@\spaceskip\z@\xspaceskip\z@
 \leavevmode{#1}\footstrut\ignorespaces#2\unskip\lower\dp\strutbox
 \vbox to\dp\strutbox{}}}
\newcount\footmarkcount@
\footmarkcount@\z@
\def\footnotemark{\let\@sf\empty\relaxnext@
 \ifhmode\edef\@sf{\spacefactor\the\spacefactor}\/\fi
 \DN@{\ifx[\next\let\next@\nextii@\else
  \ifx"\next\let\next@\nextiii@\else
  \let\next@\nextiv@\fi\fi\next@}%
 \DNii@[##1]{\footmarkform@{##1}\@sf}%
 \def\nextiii@"##1"{{##1}\@sf}%
 \def\nextiv@{\iffirstchoice@\global\advance\footmarkcount@\@ne\fi
  \footmarkform@{\number\footmarkcount@}\@sf}%
 \FN@\next@}
\def\footnotetext{\relaxnext@
 \DN@{\ifx[\next\let\next@\nextii@\else
  \ifx"\next\let\next@\nextiii@\else
  \let\next@\nextiv@\fi\fi\next@}%
 \DNii@[##1]##2{\makefootnote@{\footmarkform@{##1}}{##2}}%
 \def\nextiii@"##1"##2{\makefootnote@{##1}{##2}}%
 \def\nextiv@##1{\makefootnote@{\footmarkform@{\number\footmarkcount@}}{##1}}%
 \FN@\next@}
\def\footnote{\let\@sf\empty\relaxnext@
 \ifhmode\edef\@sf{\spacefactor\the\spacefactor}\/\fi
 \DN@{\ifx[\next\let\next@\nextii@\else
  \ifx"\next\let\next@\nextiii@\else
  \let\next@\nextiv@\fi\fi\next@}%
 \DNii@[##1]##2{\footnotemark[##1]\footnotetext[##1]{##2}}%
 \def\nextiii@"##1"##2{\footnotemark"##1"\footnotetext"##1"{##2}}%
 \def\nextiv@##1{\footnotemark\footnotetext{##1}}%
 \FN@\next@}
\def\adjustfootnotemark#1{\advance\footmarkcount@#1\relax}
\def\footnoterule{\kern-3\p@
  \hrule width 5pc\kern 2.6\p@} 
\def\captionfont@{\smc}
\def\topcaption#1#2\endcaption{%
  {\dimen@\hsize \advance\dimen@-\captionwidth@
   \rm\raggedcenter@ \advance\leftskip.5\dimen@ \rightskip\leftskip
  {\captionfont@#1}%
  \if\notempty{#2}.\enspace\ignorespaces#2\fi
  \endgraf}\nobreak\bigskip}
\def\botcaption#1#2\endcaption{%
  \nobreak\bigskip
  \setboxz@h{\captionfont@#1\if\notempty{#2}.\enspace\rm#2\fi}%
  {\dimen@\hsize \advance\dimen@-\captionwidth@
   \leftskip.5\dimen@ \rightskip\leftskip
   \noindent \ifdim\wdz@>\captionwidth@ 
   \else\hfil\fi 
  {\captionfont@#1}\if\notempty{#2}.\enspace\rm#2\fi\endgraf}}
\def\@ins{\par\begingroup\def\vspace##1{\vskip##1\relax}%
  \def\captionwidth##1{\captionwidth@##1\relax}%
  \setbox\z@\vbox\bgroup} % start a \vbox
\def\block{\RIfMIfI@\nondmatherr@\block\fi
       \else\ifvmode\vskip\abovedisplayskip\noindent\fi
        $$\def\endblock{\par\egroup$$}\fi
  \vbox\bgroup\advance\hsize-2\indenti\noindent}
\def\endblock{\par\egroup}
\def\cite#1{{\rm[{\citefont@\m@th#1}]}}
\def\citefont@{\rm}
\def\refsfont@{\eightpoint}
\outer\def\Refs{\runaway@{proclaim}%
 \relaxnext@ \DN@{\ifx\next\nofrills\DN@\nofrills{\nextii@}\else
  \DN@{\nextii@{References}}\fi\next@}%
 \DNii@##1{\penaltyandskip@{-200}\aboveheadskip
  \line{\hfil\headfont@\ignorespaces##1\unskip\hfil}\nobreak
  \vskip\belowheadskip
  \begingroup\refsfont@\sfcode`.=\@m}%
 \FN@\next@}
\def\endRefs{\par\endgroup}
\newbox\nobox@            \newbox\keybox@           \newbox\bybox@
\newbox\paperbox@         \newbox\paperinfobox@     \newbox\jourbox@
\newbox\volbox@           \newbox\issuebox@         \newbox\yrbox@
\newbox\pagesbox@         \newbox\bookbox@          \newbox\bookinfobox@
\newbox\publbox@          \newbox\publaddrbox@      \newbox\finalinfobox@
\newbox\edsbox@           \newbox\langbox@
\newif\iffirstref@        \newif\iflastref@
\newif\ifprevjour@        \newif\ifbook@            \newif\ifprevinbook@
\newif\ifquotes@          \newif\ifbookquotes@      \newif\ifpaperquotes@
\newdimen\bysamerulewd@
\setboxz@h{\refsfont@\kern3em}
\bysamerulewd@\wdz@
\newdimen\refindentwd
\setboxz@h{\refsfont@ 00. }
\refindentwd\wdz@
\outer\def\ref{\begingroup \noindent\hangindent\refindentwd
 \firstref@true \def\nofrills{\def\refkern@{\kern3sp}}%
 \ref@}
\def\ref@{\book@false \bgroup\let\endrefitem@\egroup \ignorespaces}
\def\moreref{\endrefitem@\endref@\firstref@false\ref@}%
\def\transl{\endrefitem@\endref@\firstref@false
  \book@false
  \prepunct@
  \setboxz@h\bgroup \aftergroup\unhbox\aftergroup\z@
    \def\endrefitem@{\unskip\refkern@\egroup}\ignorespaces}%
\def\emptyifempty@{\dimen@\wd\currbox@
  \advance\dimen@-\wd\z@ \advance\dimen@-.1\p@
  \ifdim\dimen@<\z@ \setbox\currbox@\copy\voidb@x \fi}
\let\refkern@\relax
\def\endrefitem@{\unskip\refkern@\egroup
  \setboxz@h{\refkern@}\emptyifempty@}\ignorespaces
\def\refdef@#1#2#3{\edef\next@{\noexpand\endrefitem@
  \let\noexpand\currbox@\csname\expandafter\eat@\string#1box@\endcsname
    \noexpand\setbox\noexpand\currbox@\hbox\bgroup}%
  \toks@\expandafter{\next@}%
  \if\notempty{#2#3}\toks@\expandafter{\the\toks@
  \def\endrefitem@{\unskip#3\refkern@\egroup
  \setboxz@h{#2#3\refkern@}\emptyifempty@}#2}\fi
  \toks@\expandafter{\the\toks@\ignorespaces}%
  \edef#1{\the\toks@}}
\refdef@\no{}{. }
\refdef@\key{[\m@th}{] }
\refdef@\by{}{}
\def\bysame{\by\hbox to\bysamerulewd@{\hrulefill}\thinspace
   \kern0sp}
\def\manyby{\message{\string\manyby is no longer necessary; \string\by
  can be used instead, starting with version 2.0 of \styname.STY}\by}
\refdef@\paper{\ifpaperquotes@``\fi\it}{}
\refdef@\paperinfo{}{}
\def\jour{\endrefitem@\let\currbox@\jourbox@
  \setbox\currbox@\hbox\bgroup
  \def\endrefitem@{\unskip\refkern@\egroup
    \setboxz@h{\refkern@}\emptyifempty@
    \ifvoid\jourbox@\else\prevjour@true\fi}%
\ignorespaces}
\refdef@\vol{\ifbook@\else\bf\fi}{}
\refdef@\issue{no. }{}
\refdef@\yr{}{}
\refdef@\pages{}{}
\def\page{\endrefitem@\def\pp@{\def\pp@{pp.~}p.~}\let\currbox@\pagesbox@
  \setbox\currbox@\hbox\bgroup\ignorespaces}
\def\pp@{pp.~}
\def\book{\endrefitem@ \let\currbox@\bookbox@
 \setbox\currbox@\hbox\bgroup\def\endrefitem@{\unskip\refkern@\egroup
  \setboxz@h{\ifbookquotes@``\fi}\emptyifempty@
  \ifvoid\bookbox@\else\book@true\fi}%
  \ifbookquotes@``\fi\it\ignorespaces}
\def\inbook{\endrefitem@
  \let\currbox@\bookbox@\setbox\currbox@\hbox\bgroup
  \def\endrefitem@{\unskip\refkern@\egroup
  \setboxz@h{\ifbookquotes@``\fi}\emptyifempty@
  \ifvoid\bookbox@\else\book@true\previnbook@true\fi}%
  \ifbookquotes@``\fi\ignorespaces}
\refdef@\eds{(}{, eds.)}
\def\ed{\endrefitem@\let\currbox@\edsbox@
 \setbox\currbox@\hbox\bgroup
 \def\endrefitem@{\unskip, ed.)\refkern@\egroup
  \setboxz@h{(, ed.)}\emptyifempty@}(\ignorespaces}
\refdef@\bookinfo{}{}
\refdef@\publ{}{}
\refdef@\publaddr{}{}
\refdef@\finalinfo{}{}
\refdef@\lang{(}{)}

\let\refdef@\relax 
\def\ppunbox@#1{\ifvoid#1\else\prepunct@\unhbox#1\fi}
\def\nocomma@#1{\ifvoid#1\else\changepunct@3\prepunct@\unhbox#1\fi}
\def\changepunct@#1{\ifnum\lastkern<3 \unkern\kern#1sp\fi}
\def\prepunct@{\count@\lastkern\unkern
  \ifnum\lastpenalty=0
    \let\penalty@\relax
  \else
    \edef\penalty@{\penalty\the\lastpenalty\relax}%
  \fi
  \unpenalty
  \let\refspace@\ \ifcase\count@,% usual case, do a comma
\or;\or.\or % do nothing; this case is from nofrills.
  \or\let\refspace@\relax
  \else,\fi
  \ifquotes@''\quotes@false\fi \penalty@ \refspace@
}
\def\transferpenalty@#1{\dimen@\lastkern\unkern
  \ifnum\lastpenalty=0\unpenalty\let\penalty@\relax
  \else\edef\penalty@{\penalty\the\lastpenalty\relax}\unpenalty\fi
  #1\penalty@\kern\dimen@}
\def\endref{\endrefitem@\lastref@true\endref@
  \par\endgroup \prevjour@false \previnbook@false }
\def\endref@{%
\iffirstref@
  \ifvoid\nobox@\ifvoid\keybox@\indent\fi
  \else\hbox to\refindentwd{\hss\unhbox\nobox@}\fi
  \ifvoid\keybox@
  \else\ifdim\wd\keybox@>\refindentwd
         \box\keybox@
       \else\hbox to\refindentwd{\unhbox\keybox@\hfil}\fi\fi
  \kern4sp\ppunbox@\bybox@
\fi 
  \ifvoid\paperbox@
  \else\prepunct@\unhbox\paperbox@
    \ifpaperquotes@\quotes@true\fi\fi
  \ppunbox@\paperinfobox@
  \ifvoid\jourbox@
    \ifprevjour@ \nocomma@\volbox@
      \nocomma@\issuebox@
      \ifvoid\yrbox@\else\changepunct@3\prepunct@(\unhbox\yrbox@
        \transferpenalty@)\fi
      \ppunbox@\pagesbox@
    \fi 
  \else \prepunct@\unhbox\jourbox@
    \nocomma@\volbox@
    \nocomma@\issuebox@
    \ifvoid\yrbox@\else\changepunct@3\prepunct@(\unhbox\yrbox@
      \transferpenalty@)\fi
    \ppunbox@\pagesbox@
  \fi 
  \ifbook@\prepunct@\unhbox\bookbox@ \ifbookquotes@\quotes@true\fi \fi
  \nocomma@\edsbox@
  \ppunbox@\bookinfobox@
  \ifbook@\ifvoid\volbox@\else\prepunct@ vol.~\unhbox\volbox@
  \fi\fi
  \ppunbox@\publbox@ \ppunbox@\publaddrbox@
  \ifbook@ \ppunbox@\yrbox@
    \ifvoid\pagesbox@
    \else\prepunct@\pp@\unhbox\pagesbox@\fi
  \else
    \ifprevinbook@ \ppunbox@\yrbox@
      \ifvoid\pagesbox@\else\prepunct@\pp@\unhbox\pagesbox@\fi
    \fi \fi
  \ppunbox@\finalinfobox@
  \iflastref@
    \ifvoid\langbox@.\ifquotes@''\fi
    \else\changepunct@2\prepunct@\unhbox\langbox@\fi
  \else
    \ifvoid\langbox@\changepunct@1%
    \else\changepunct@3\prepunct@\unhbox\langbox@
      \changepunct@1\fi
  \fi
}
\outer\def\enddocument{%
 \runaway@{proclaim}%
\ifmonograph@ % do nothing
\else
 \nobreak
 \thetranslator@
 \count@\z@ \loop\ifnum\count@<\addresscount@\advance\count@\@ne
 \csname address\number\count@\endcsname
 \csname email\number\count@\endcsname
 \repeat
\fi
 \vfill\supereject\end}

%Modifizierte Fonts fr Proclaim,...
\def\headfont@{\headfonts}
\def\proclaimheadfont@{\bf}
\def\specialheadfont@{\bf}
\def\subheadfont@{\bf}
\def\demoheadfont@{\smc}

%Kontrollsequenzen fr Inhaltsverzeichnis und Index:
\newif\ifThisToToc \ThisToTocfalse
\newif\iftocloaded \tocloadedfalse

\def\C@L{\noexpand\Cal}\def\B@B{\noexpand\Bbb}\def\fR@K{\noexpand\frak}
\def\S@{\noexpand\S}\def\P@P{\noexpand\"}
\def\xpar{\\}

\def\writetoc#1{\iftocloaded\ifThisToToc\begingroup\def\totoc{}
  \def\Cal{\noexpand\C@L}\def\Bbb{\noexpand\B@B}
  \def\frak{\noexpand\fR@K}\def\goth{\frak}\def\S{\noexpand\S@}
  \def\"{\noexpand\P@P}
  \def\xpar{\par\penalty100000 }\def\idx##1{##1}\def\\{\xpar}
  \edef\next@{\write\toc{\noindent#1\leaderfill\noexpand\folio\par}}%
  \next@\endgroup\global\ThisToTocfalse\fi\fi}
\def\leaderfill{\leaders\hbox to 1em{\hss.\hss}\hfill}

\newif\ifindexloaded \indexloadedfalse
\def\idx#1{\ifindexloaded\begingroup\def\ign{}\def\it{}\def\/{}%
 \def\smc{}\def\bf{}\def\tt{}%
 \def\Cal{\noexpand\C@L}\def\Bbb{\noexpand\B@B}%
 \def\frak{\noexpand\fR@K}\def\goth{\frak}\def\S{\noexpand\S@}%
  \def\"{\noexpand\P@P}%
 {\edef\next@{\write\index{#1, \noexpand\folio}}\next@}%
 \endgroup\fi{#1}}
\def\ign#1{}

\def\totoc{\global\ThisToToctrue}

\outer\def\head#1\endhead{\par\penaltyandskip@{-200}\aboveheadskip
 {\headfont@\raggedcenter@\interlinepenalty\@M
 \ignorespaces#1\endgraf}\nobreak
 \vskip\belowheadskip
 \headmark{#1}\writetoc{#1}}

\outer\def\chaphead#1\endchaphead{\par\penaltyandskip@{-200}\aboveheadskip
 {\chapheadfonts\raggedcenter@\interlinepenalty\@M
 \ignorespaces#1\endgraf}\nobreak
 \vskip3\belowheadskip
 \headmark{#1}\writetoc{#1}}

\def\folio{{\foliofont@\ifnum\pageno<\z@ \romannumeral-\pageno
 \else\number\pageno \fi}}
\newtoks\leftheadtoks
\newtoks\rightheadtoks

%Uppercase ist abgestellt:
\def\leftheadtext{\nofrills@{\relax}\lht@
  \DNii@##1{\leftheadtoks\expandafter{\lht@{##1}}%
    \mark{\the\leftheadtoks\noexpand\else\the\rightheadtoks}
    \ifsyntax@\setboxz@h{\def\\{\unskip\space\ignorespaces}%
        \headlinefont@##1}\fi}%
  \FN@\next@}
%Uppercase ist abgestellt:
\def\rightheadtext{\nofrills@{\relax}\rht@
  \DNii@##1{\rightheadtoks\expandafter{\rht@{##1}}%
    \mark{\the\leftheadtoks\noexpand\else\the\rightheadtoks}%
    \ifsyntax@\setboxz@h{\def\\{\unskip\space\ignorespaces}%
        \headlinefont@##1}\fi}%
  \FN@\next@}
%\headline={\def\chapter#1{\chapterno@. }%
%  \def\\{\unskip\space\ignorespaces}\headlinefont@
%  \ifodd\pageno \rightheadline \else \leftheadline\fi}
\def\NoRunningHeads{\global\runheads@false\global\let\headmark\eat@}

\newif\iffirstpage@     \firstpage@true
\newif\ifrunheads@      \runheads@true

%Ergnzungen zu Runningheads und Pagenumbers:
\newdimen\fullhsize \fullhsize=\hsize
\newdimen\fullvsize \fullvsize=\vsize
\def\fullline{\hbox to\fullhsize}

\def\pagenumbers{\gdef\folio{\folio@}}

\let\norunningheads\NoRunningHeads
\def\userunningheads{\global\runheads@true}
%Default: Seitennumerierung unten.
\norunningheads

\headline={\def\chapter#1{\chapterno@. }%
  \def\\{\unskip\space\ignorespaces}\ifrunheads@\headlinefont@
    \ifodd\pageno\rightheadline \else\leftheadline\fi
   \else\hfil\fi\ifNoRunHeadline\global\NoRunHeadlinefalse\fi}
\let\folio@\folio
\def\foliofont@{\foliofont}
\def\foliofont{\eightrm}
\def\headlinefont@{\headlinefont}
\def\headlinefont{\eightpoint\smc}
\def\leftheadline{\rlap{\folio}\hfill
   \ifNoRunHeadline\else\iftrue\topmark\fi\fi \hfill}
\def\rightheadline{\hfill\ifNoRunHeadline
   \else \expandafter\iffalse\botmark\fi\fi
  \hfill \llap{\folio}}
\footline={{\eightpoint\bottremark}%
   \ifrunheads@\else\hfil{\let\foliofont\tenrm\folio}\fi\hfil}
\def\bottremark{}
 
%Definition von \norunninghead:
\newif\ifNoRunHeadline      
\def\norunninghead{\global\NoRunHeadlinetrue}
\norunninghead

\output={\output@}
%\def\output@{\shipout\vbox{%
% \iffirstpage@ \global\firstpage@false
%  \pagebody \logo@ \makefootline%
% \else \ifrunheads@ \makeheadline \pagebody
%       \else \pagebody \makefootline \fi
% \fi}%
% \advancepageno \ifnum\outputpenalty>-\@MM\else\dosupereject\fi}
%
%Modifizierter Output + Index-Output:
\newif\ifoffset\offsetfalse
\output={\output@}
\def\output@{%
 \ifoffset 
  \ifodd\count0\advance\hoffset by0.5truecm
   \else\advance\hoffset by-0.5truecm\fi\fi
 \shipout\vbox{%
  \makeheadline \pagebody \makefootline }%
 \advancepageno \ifnum\outputpenalty>-\@MM\else\dosupereject\fi}

\def\indexoutput#1{%
  \ifoffset 
   \ifodd\count0\advance\hoffset by0.5truecm
    \else\advance\hoffset by-0.5truecm\fi\fi
  \shipout\vbox{\makeheadline
  \vbox to\fullvsize{\boxmaxdepth\maxdepth%
  \ifvoid\topins\else\unvbox\topins\fi% 
  #1 %
  \ifvoid\footins\else % footnote info is present
    \vskip\skip\footins
    \footnoterule
    \unvbox\footins\fi
  \ifr@ggedbottom \kern-\dimen@ \vfil \fi}%
  \baselineskip2pc
  \makefootline}%
 \global\advance\pageno\@ne
 \ifnum\outputpenalty>-\@MM\else\dosupereject\fi}
 
 \newbox\partialpage \newdimen\halfsize \halfsize=0.5\fullhsize
 \advance\halfsize by-0.5em

 \def\begindoublecolumns{\output={\indexoutput{\unvbox255}}%
   \begingroup \def\line{\fullline}
   \output={\global\setbox\partialpage=\vbox{\unvbox255\bigskip}}\eject
   \output={\doublecolumnout}\hsize=\halfsize \vsize=2\fullvsize}
 \def\enddoublecolumns{\output={\balancecolumns}\eject
  \endgroup \pagegoal=\fullvsize%
  \output={\output@}}
\def\doublecolumnout{\splittopskip=\topskip \splitmaxdepth=\maxdepth
  \dimen@=\fullvsize \advance\dimen@ by-\ht\partialpage
  \setbox0=\vsplit255 to \dimen@ \setbox2=\vsplit255 to \dimen@
  \indexoutput{\pagesofar} \unvbox255 \penalty\outputpenalty}
\def\pagesofar{\unvbox\partialpage
  \wd0=\hsize \wd2=\hsize \hbox to\fullhsize{\box0\hfil\box2}}
\def\balancecolumns{\setbox0=\vbox{\unvbox255} \dimen@=\ht0
  \advance\dimen@ by\topskip \advance\dimen@ by-\baselineskip
  \divide\dimen@ by2 \splittopskip=\topskip
  {\vbadness=10000 \loop \global\setbox3=\copy0
    \global\setbox1=\vsplit3 to\dimen@
    \ifdim\ht3>\dimen@ \global\advance\dimen@ by1pt \repeat}
  \setbox0=\vbox to\dimen@{\unvbox1} \setbox2=\vbox to\dimen@{\unvbox3}
  \pagesofar}

\tenpoint
\catcode`\@=\active

\def\smallheadings{\let\chapheadfonts\tenpoint\let\headfonts\tenpoint}

\tenpoint
\catcode`\@=\active

\def\LL{\leavevmode\setbox0=\hbox{L}\hbox to\wd0{\hss\char'40L}}
\def\al{\alpha}
\def\be{\beta}

\def\de{\delta}
\def\ep{\varepsilon}

\def\la{\lambda}

\def\si{\sigma}
\def\ta{\tau}

\def\om{\omega}

\def\La{\Lambda}

            %used for crossreferencing, Tex should ignore.
             %used for refencing (section-numbers)
          %used for new-section numbers

\def\today{\ifcase\month\or
 January\or February\or March\or April\or May\or June\or
 July\or August\or September\or October\or November\or December\fi
 \space\number\day, \number\year}
 %zum Nummerieren

\def\({\left(}
\def\){\right)}
\def\[{\left[}
\def\]{\right]}

\def\sgn{\operatorname{sgn}}

\def\3{\ss}
\catcode`\@=11
\def\dddot#1{\vbox{\ialign{##\crcr
      .\hskip-.5pt.\hskip-.5pt.\crcr\noalign{\kern1.5\p@\nointerlineskip}
      $\hfil\displaystyle{#1}\hfil$\crcr}}}

\newif\iftab@\tab@false
\newif\ifvtab@\vtab@false
\def\tab{\bgroup\tab@true\vtab@false\vst@bfalse\Strich@false%
   \def\\{\global\hline@@false%
     \ifhline@\global\hline@false\global\hline@@true\fi\cr}
   \edef\l@{\the\leftskip}\ialign\bgroup\hskip\l@##\hfil&&##\hfil\cr}
\def\endtab{\cr\egroup\egroup}
\def\vtab{\vtop\bgroup\vst@bfalse\vtab@true\tab@true\Strich@false%
   \bgroup\def\\{\cr}\ialign\bgroup&##\hfil\cr}
\def\endvtab{\cr\egroup\egroup\egroup}
\def\stab{\D@cke0.5pt\null 
 \bgroup\tab@true\vtab@false\vst@bfalse\Strich@true\Let@@\vspace@
 \normalbaselines\offinterlineskip
  \openup\spreadmlines@
 \edef\l@{\the\leftskip}\ialign
 \bgroup\hskip\l@##\hfil&&##\hfil\crcr}
\def\endstab{\crcr\egroup
 \egroup}
\newif\ifvst@b\vst@bfalse
\def\vstab{\D@cke0.5pt\null
 \vtop\bgroup\tab@true\vtab@false\vst@btrue\Strich@true\bgroup\Let@@\vspace@
 \normalbaselines\offinterlineskip
  \openup\spreadmlines@\bgroup}
\def\endvstab{\crcr\egroup\egroup
 \egroup\tab@false\Strich@false}

\newdimen\htstrut@
\htstrut@8.5\p@
\newdimen\htStrut@
\htStrut@12\p@
\newdimen\dpstrut@
\dpstrut@3.5\p@
\newdimen\dpStrut@
\dpStrut@3.5\p@
\def\openup{\afterassignment\@penup\dimen@=}
\def\@penup{\advance\lineskip\dimen@
  \advance\baselineskip\dimen@
  \advance\lineskiplimit\dimen@
  \divide\dimen@ by2
  \advance\htstrut@\dimen@
  \advance\htStrut@\dimen@
  \advance\dpstrut@\dimen@
  \advance\dpStrut@\dimen@}
\def\Let@@{\relax\iffalse{\fi%
    \def\\{\global\hline@@false%
     \ifhline@\global\hline@false\global\hline@@true\fi\cr}%
    \iffalse}\fi}
\def\matrix{\null\,\vcenter\bgroup
 \tab@false\vtab@false\vst@bfalse\Strich@false\Let@@\vspace@
 \normalbaselines\openup\spreadmlines@\ialign
 \bgroup\hfil$\m@th##$\hfil&&\quad\hfil$\m@th##$\hfil\crcr
 \Mathstrut@\crcr\noalign{\kern-\baselineskip}}
\def\endmatrix{\crcr\Mathstrut@\crcr\noalign{\kern-\baselineskip}\egroup
 \egroup\,}
\def\smatrix{\D@cke0.5pt\null\,
 \vcenter\bgroup\tab@false\vtab@false\vst@bfalse\Strich@true\Let@@\vspace@
 \normalbaselines\offinterlineskip
  \openup\spreadmlines@\ialign
 \bgroup\hfil$\m@th##$\hfil&&\quad\hfil$\m@th##$\hfil\crcr}
\def\endsmatrix{\crcr\egroup
 \egroup\,\Strich@false}
\newdimen\D@cke
\def\Dicke#1{\global\D@cke#1}
\newtoks\tabs@\tabs@{&}
\newif\ifStrich@\Strich@false
\newif\iff@rst

\def\Stricherr@{\iftab@\ifvtab@\errmessage{\noexpand\s not allowed
     here. Use \noexpand\vstab!}%
  \else\errmessage{\noexpand\s not allowed here. Use \noexpand\stab!}%
  \fi\else\errmessage{\noexpand\s not allowed
     here. Use \noexpand\smatrix!}\fi}
\def\format{\ifvst@b\else\crcr\fi\egroup\iffalse{\fi\ifnum`}=0 \fi\format@}
\def\format@#1\\{\def\preamble@{#1}%
 \def\Str@chfehlt##1{\ifx##1\s\Stricherr@\fi\ifx##1\\\let\Next\relax%
   \else\let\Next\Str@chfehlt\fi\Next}%
 \def\c{\hfil\noexpand\ifhline@@\hbox{\vrule height\htStrut@%
   depth\dpstrut@ width\z@}\noexpand\fi%
   \ifStrich@\hbox{\vrule height\htstrut@ depth\dpstrut@ width\z@}%
   \fi\iftab@\else$\m@th\fi\the\hashtoks@\iftab@\else$\fi\hfil}%
 \def\r{\hfil\noexpand\ifhline@@\hbox{\vrule height\htStrut@%
   depth\dpstrut@ width\z@}\noexpand\fi%
   \ifStrich@\hbox{\vrule height\htstrut@ depth\dpstrut@ width\z@}%
   \fi\iftab@\else$\m@th\fi\the\hashtoks@\iftab@\else$\fi}%
 \def\l{\noexpand\ifhline@@\hbox{\vrule height\htStrut@%
   depth\dpstrut@ width\z@}\noexpand\fi%
   \ifStrich@\hbox{\vrule height\htstrut@ depth\dpstrut@ width\z@}%
   \fi\iftab@\else$\m@th\fi\the\hashtoks@\iftab@\else$\fi\hfil}%
 \def\s{\ifStrich@\ \the\tabs@\vrule width\D@cke\the\hashtoks@%
          \fi\the\tabs@\ }%
 \def\sa{\ifStrich@\vrule width\D@cke\the\hashtoks@%
            \the\tabs@\ %
            \fi}%
 \def\se{\ifStrich@\ \the\tabs@\vrule width\D@cke\the\hashtoks@\fi}%
 \def\cd{\hfil\noexpand\ifhline@@\hbox{\vrule height\htStrut@%
   depth\dpstrut@ width\z@}\noexpand\fi%
   \ifStrich@\hbox{\vrule height\htstrut@ depth\dpstrut@ width\z@}%
   \fi$\dsize\m@th\the\hashtoks@$\hfil}%
 \def\rd{\hfil\noexpand\ifhline@@\hbox{\vrule height\htStrut@%
   depth\dpstrut@ width\z@}\noexpand\fi%
   \ifStrich@\hbox{\vrule height\htstrut@ depth\dpstrut@ width\z@}%
   \fi$\dsize\m@th\the\hashtoks@$}%
 \def\ld{\noexpand\ifhline@@\hbox{\vrule height\htStrut@%
   depth\dpstrut@ width\z@}\noexpand\fi%
   \ifStrich@\hbox{\vrule height\htstrut@ depth\dpstrut@ width\z@}%
   \fi$\dsize\m@th\the\hashtoks@$\hfil}%
 \ifStrich@\else\Str@chfehlt#1\\\fi%
 \setbox\z@\hbox{\xdef\Preamble@{\preamble@}}\ifnum`{=0 \fi\iffalse}\fi
 \ialign\bgroup\span\Preamble@\crcr}
\newif\ifhline@\hline@false
\newif\ifhline@@\hline@@false
\def\hlinefor#1{\multispan@{\strip@#1 }\leaders\hrule height\D@cke\hfill%
    \global\hline@true\ignorespaces}
\def\Item "#1"{\par\noindent\hangindent2\parindent%
  \hangafter1\setbox0\hbox{\rm#1\enspace}\ifdim\wd0>2\parindent%
  \box0\else\hbox to 2\parindent{\rm#1\hfil}\fi\ignorespaces}
\def\ITEM #1"#2"{\par\noindent\hangafter1\hangindent#1%
  \setbox0\hbox{\rm#2\enspace}\ifdim\wd0>#1%
  \box0\else\hbox to 0pt{\rm#2\hss}\hskip#1\fi\ignorespaces}
\def\item"#1"{\par\noindent\hang%
  \setbox0=\hbox{\rm#1\enspace}\ifdim\wd0>\the\parindent%
  \box0\else\hbox to \parindent{\rm#1\hfil}\enspace\fi\ignorespaces}
\let\plainitem@\item
\catcode`\@=13

\hsize13cm
\vsize19cm
\newdimen\fullhsize
\newdimen\fullvsize
\newdimen\halfsize
\fullhsize13cm
\fullvsize19cm
\halfsize=0.5\fullhsize
\advance\halfsize by-0.5em

\magnification1200

\TagsOnRight

\def\BaDSAA{1}
\def\BrHiAA{2}
\def\ChriAA{3}
\def\CJKrAB{4}
\def\DoFSAA{5}
\def\ElouAA{6}
\def\IsmaAA{7}
\def\KratBN{8}
\def\KratCI{9}
\def\KratCP{10}
\def\LascAZ{11}
\def\MuWYAA{12}
\def\StanBI{13}
\def\SzegAA{14}
\def\VienAE{15}

%equation numbers
\def\AA{1.1}
\def\AB{1.2}
\def\ABa{1.3}
\def\AC{1.4}

\def\AZ{2.1}
\def\AZa{2.2}
\def\AZb{2.3}
\def\AZc{2.4}
\def\BA{3.1}
\def\BB{3.2}
\def\BC{3.3}
\def\BD{3.4}
\def\BE{3.5}
\def\BF{3.6}
\def\BG{3.7}
\def\BH{3.8}

\def\CA{4.1}
\def\CAa{4.2}
\def\CB{4.3}
\def\CC{4.4}
\def\CD{4.5}

\def\DAd{5.1}
\def\DAa{5.2}
\def\DAb{5.3}
\def\DAc{5.4}
\def\DA{5.5}
\def\DB{5.6}
\def\DC{5.7}

\def\YA{6.1}
\def\YB{6.2}
\def\YC{6.3}
\def\YD{6.4}

\def\EA{7.1}
\def\EAa{7.2}
\def\EB{7.3}
\def\EC{7.4}

\def\ZA{8.1}
\def\ZB{8.2}
\def\ZBa{8.3}
\def\ZC{8.4}
\def\ZD{8.5}
\def\ZE{8.6}
\def\ZF{8.7}
\def\ZG{8.8}
\def\ZH{8.9}

%

%theorem numbers
\def\TA{1}
\def\TB{2}
\def\TC{3}
\def\TD{4}
\def\TEa{5}
\def\TE{6}
\def\TF{7}
\def\TFa{8}
\def\TG{9}

%

%figure numbers

\def\LHS{\operatorname{LHS}}

\topmatter 
\title Hankel determinants of linear combinations
of moments of orthogonal polynomials, II
\endtitle 
\author C.~Krattenthaler
\endauthor 
\affil 
Fakult\"at f\"ur Mathematik, Universit\"at Wien,\\
Oskar-Morgenstern-Platz~1, A-1090 Vienna, Austria.\\
WWW: \tt http://www.mat.univie.ac.at/\~{}kratt
\endaffil
\address Fakult\"at f\"ur Mathematik, Universit\"at Wien,
Oskar-Morgenstern-Platz~1, A-1090 Vienna, Austria.\newline
http://www.mat.univie.ac.at/\~{}kratt
\endaddress

\dedicatory
Dedicated to the memory of Richard Askey
\enddedicatory

\thanks Research partially supported by the Austrian
Science Foundation FWF (grant S50-N15)
in the framework of the Special Research Program
``Algorithmic and Enumerative Combinatorics"%
\endthanks

\subjclass Primary 33C45;
 Secondary 05A15 11C20 15A15 42C05
\endsubjclass
\keywords Hankel determinants, moments of orthogonal polynomials,
Catalan numbers, 
Motzkin numbers, 
Schr\"oder numbers, Riordan numbers, Fine numbers,
central binomial coefficients, central trinomial numbers,
Delannoy numbers, Chebyshev polynomials, Dodgson condensation
\endkeywords
\abstract 
We present a formula that expresses the Hankel determinants of
a linear combination of length~$d+1$ of moments of orthogonal polynomials
in terms of a $d\times d$ determinant of the orthogonal polynomials.
This formula exists somehow hidden in the folklore of the theory of
orthogonal polynomials but deserves to be better known, and be
presented correctly and with full proof. We present four
fundamentally different
proofs, one that uses classical formulae from the theory of orthogonal
polynomials, one that uses a vanishing argument and is due to
Elouafi [{\it J. Math\. Anal\. Appl\.} {\bf 431} (2015), 1253--1274]
(but given in an incomplete form there), one that is inspired by
random matrix theory and is due to Br\'ezin and Hikami
[{\it Comm\. Math\. Phys\.} {\bf214} (2000), 111--135],
and one that uses (Dodgson) condensation.
We give two applications of the formula. In the first application,
we explain how to
compute such Hankel determinants in a singular case. The second
application concerns the linear recurrence of such Hankel
determinants for a certain class of moments that covers numerous
classical combinatorial sequences, including Catalan numbers,
Motzkin numbers,
central binomial coefficients, central trinomial coefficients,
central Delannoy numbers, Schr\"oder numbers, Riordan numbers,
and Fine numbers.
\endabstract
\endtopmatter
\document

\subhead 1. Introduction\endsubhead
The purpose of this article is to put to the fore a fundamental
formula for orthogonal polynomials that 
is implicitly hidden in the classical literature on
orthogonal polynomials.
It is so well hidden that seemingly
even top experts of the theory of orthogonal
polynomials are not aware of the formula.
How and why this is possible is explained in greater detail in
Section~2. My literature search led me to discover that the formula
is stated in Lascoux's book \cite{\LascAZ}, albeit incorrectly, but
with a correct proof. Subsequently, I realised that the formula is
stated correctly by Elouafi in \cite{\ElouAA}, however with an
incomplete proof. Finally, in reaction to~\cite{\KratCP},
Arno Kuijlaars pointed out to me that the formula appears in~\cite{\BaDSAA},
in which a result due to Br\'ezin and Hikami \cite{\BrHiAA} is cited. 
Both papers
contain correct statement and (different) proofs, however
they use random matrix language.
Again, see Section~2 for more details.

So, let me present this formula without further ado.
Let $\big(p_n(x)\big)_{n\ge0}$ be a sequence of monic polynomials
over a field~$K$ of characteristic zero\footnote{For the analyst, 
(usually) this field is the field of 
real numbers, and a further restriction is that the linear
functional~$L$ is defined by a measure with non-negative density. 
However, the formulae
in this paper do not need these restrictions and 
are valid in this wider context of ``formal orthogonality".} with
$\deg p_n(x)=n$, and assume that they are orthogonal with respect to the linear
functional~$L$, i.e., they satisfy $L(p_m(x)p_n(x))=\om_n\de_{m,n}$ with
$\om_n\ne0$ for all~$n$,
where $\de_{m,n}$ is the Kronecker delta.
Furthermore, we write $\mu_n$ for the $n$-th moment
$L(x^n)$ of the functional~$L$, for which we also use the {\it umbral notation}
$\mu^n\equiv \mu_n$.\footnote{``Umbral notation" means that an
expression that is a polynomial in~$\mu$ is expanded out, and then
every occurrence of $\mu^n$ is replaced by~$\mu_n$. So, for example,
the umbral expression $\mu^2(x_1+\mu)(x_2+\mu)$ means
$x_1x_2\mu_2+(x_1+x_2)\mu_3+\mu_4$.} 

\proclaim{Theorem \TA}
Let $n$ and $d$ be non-negative integers.
Given variables $x_1,x_2,\dots,x_d$, and using the above explained
umbral notation, we have
$$
\frac {\det\limits_{0\le i,j\le n-1}\left(
\mu^{i+j}\prod _{\ell=1} ^{d}(x_\ell+\mu)\right)}
{\det\limits_{0\le i,j\le n-1}\left(\mu_{i+j}\right)}
=(-1)^{nd}
\frac {\displaystyle\det_{1\le i,j\le d}\left(p_{n+i-1}(-x_j)\right)}
{
\prod\limits _{1\le i<j\le d} ^{}(x_i-x_j)}.
\tag\AA
$$
Here, determinants of empty matrices and empty products are understood
to equal~$1$.
\endproclaim

\remark{Remark}
The theory of orthogonal polynomials guarantees that in our setting
(namely due to the condition $\om_n\ne0$ in the orthogonality)
the Hankel determinant of moments in the denominator on the left-hand
side of~(\AA) is non-zero.
\endremark

We may rewrite (\AA) using quantities that appear in the
{\it three-term recurrence}
$$p_n(x)=(x-s_{n-1})p_{n-1}(x)-t_{n-2}p_{n-2}(x),
\quad \text{for }n\ge1,
\tag\AB
$$
with initial values
$p_{-1} (x)= 0$ and $p_0 (x)=1$, that is satisfied by the polynomials
according to Favard's theorem (see e.g\. \cite{\KratBN,
Theorems~11--13})
for some sequences $(s_n)_{n\ge0}$ and $(t_n)_{n\ge0}$
of elements of~$K$ with $t_n\ne0$ for all $n$.
Namely, using the well-known fact (see e.g\. \cite{\VienAE,
  Ch.~IV, Cor.~6}) 
$$
\det_{0\le i,j\le n-1}\left(\mu_{i+j}\right)
=\prod _{i=0} ^{n-1}t_i^{n-i-1},
\tag\ABa
$$
the formula (\AA) becomes
$$
\det\limits_{0\le i,j\le n-1}\left(
\mu^{i+j}\prod _{\ell=1} ^{d}(x_\ell+\mu)\right)
=(-1)^{nd+\binom d2}
\bigg(\prod _{i=0} ^{n-1}t_i^{n-i-1}\bigg)
\frac {\displaystyle\det_{1\le i,j\le d}\left(p_{n+i-1}(-x_j)\right)}
{
\prod\limits _{1\le i<j\le d} ^{}(x_j-x_i)}.
\tag\AC
$$
This form reveals that we may regard the formula as a polynomial formula in
the $x_i$'s {\it and\/} the $s_i$'s and $t_i$'s. Indeed, the determinant
on the right-hand side, being a skew-symmetric polynomial in
the~$x_i$'s, is divisible by the Vandermonde product in the denominator.

\medskip
In the next section, I will present the history of Theorem~\TA, 
from a strongly biased (personal) view. As I explain there, I
discovered the formula on my own while thinking about
Conjecture~8 in \cite{\CJKrAB}, and also came up with a proof,
presented here in Section~3. Later I found the earlier mentioned
occurrences of the formula in \cite{\LascAZ}, \cite{\ElouAA},
\cite{\BaDSAA} and~\cite{\BrHiAA}.
Lascoux's argument (the one in~\cite{\BaDSAA} is essentially the same), 
which follows the classical literature of
orthogonal polynomials (but is presented in \cite{\LascAZ} in his very personal
language), is presented in Section~4 (in ``standard" language). 
Section~5 brings the
completion of Elouafi's vanishing argument. The
random matrix-inspired proof due to Br\'ezin and Hikami is the subject
of Section~6.

Sections~7 and~8 address
issues that come from my initial motivation
(and Elouafi's) that in the end led
to the discovery of Theorem~\TA: Hankel determinants of
linear combinations of combinatorial sequences.
Section~7 addresses the case in which in (\AA) the~$x_i$'s
are all equal to each other. In that case, it is the limit
formula in Proposition~\TEa\ in Section~5 that has to be applied. 
We show in Section~7 that Elouafi's recurrence approach for that
case can be replaced by an approach yielding completely explicit
expressions. Finally, in Section~8 we show that the theory of
linear recurrent sequences with constant coefficients implies that,
in the case where the coefficients $s_i$ and $t_i$ in the three-term
recurrence~(\AB) are constant for large~$i$, the scaled Hankel
determinants of linear combinations of moments on the left-hand side
of~(\AA) satisfy a linear recurrence with constant coefficients of
order~$2^d$, plus some more specific assertions about the coefficients
in this linear recurrence, see Corollary~\TG. This proves conjectures
from \cite{\DoFSAA}, vastly generalising them.

%What is the significance of Theorem~\TA?
%
%fundamental identity in the theory of orthogonal polynomials
%-> Christoffel/Szeg\H o, Leclerc, Dur\'an
%
%Here: Hankel determinants of linear combinations of
%Motzkin-related/combinatorial numbers
%
%Anwendungen:
%* Beweis von Conjecture 8
%* C-recursiveness der Hankeldeterminanten wenn s_i's und t_i's
%  konstant sind.
%* explizite Formeln fuer Hankeldeterminanten von Linearkombinationen
%  von Motzkinzahlen, Catalanzahlen, usw.

\subhead 2. History of Theorem \TA\ --- a (very) personal view\endsubhead
I discovered Theorem~{\TA} on my own, in a very roundabout way.
It started with an email of Johann Cigler in which he asked me
for a proof of a special case of
$$
\frac {\det\left(x_1x_2 \mu_{i+j}+(x_1+x_2)\mu_{i+j+1}+\mu_{i+j+2}\right)_{i,j=0}^{n-1}} 
{\det\left(\mu_{i+j}\right)_{i,j=0}^{n-1}}
=\sum_{j=0}^np_{j}(-x_2)p_j(-x_2)
\prod _{\ell=j} ^{n-1}t_\ell.
\tag\AZ
$$
We quickly realised that we can actually prove the above identity,
which became the first main result in \cite{\CJKrAB} (see
Theorem~1 there; 
the reader should notice that the left-hand side of (\AZ) agrees with
the left-hand side of (\AA), while the right-hand sides do not agree;
in retrospect, the equality of the right-hand sides is equivalent 
to the Christoffel--Darboux identity, cf\. \cite{\SzegAA, Theorem~3.2.2}).
We then proceeded to derive a (more complicated) triple-sum expression
for the ``next" case (see \cite{\CJKrAB, Theorem~5}),
$$
\frac {\det\left(x_1x_2x_3 \mu_{i+j}+(x_1x_2+x_2x_3+x_3x_1)\mu_{i+j+1}
+(x_1+x_2+x_3)\mu_{i+j+2}+\mu_{i+j+3}\right)_{i,j=0}^{n-1}} 
{\det\left(\mu_{i+j}\right)_{i,j=0}^{n-1}}.
$$
In the special case where the $s_i$'s and the $t_i$'s are constant
for $i\ge1$,
the orthogonal polynomial $p_n(x)$ can be expressed in terms of
a linear combination of Chebyshev polynomials (see \cite{\CJKrAB,
Eq.~(4.2)}). This allowed us to
evaluate the sum on the right-hand side of (\AZ) and the
afore-mentioned triple sum. We recognised a pattern, and this led
us to conjecture a precise formula for
$$
\frac {\det\limits_{0\le i,j\le n-1}\left(
\mu^{i+j}\prod _{\ell=1} ^{d}(x_\ell+\mu)\right)}
{\det\limits_{0\le i,j\le n-1}\left(\mu_{i+j}\right)}
$$
(see \cite{\CJKrAB, Conj.~8}), again expressed in terms of Chebyshev
polynomials. Subsequently, I realised that this conjectural expression
could be simplified (by means of \cite{\CJKrAB, Eq.~(4.2)}). The
result was the right-hand side of (\AA), in the special case where
the $s_i$'s and $t_i$'s are constant for $i\ge1$. The obvious question
at that point then was: does Formula~(\AA) also hold if the $s_i$'s
and $t_i$'s are generic? Computer experiments said ``yes". 

At this point I told myself: this
identity, being a completely general identity of fundamental nature 
connecting orthogonal
polynomials and their moments, must be known. Naturally, I consulted
standard books on orthogonal polynomials, such as Szeg\H o's classic
\cite{\SzegAA},
but I could not find it. After a while I then started to think about a
proof. 
I figured out the proof of~(\AA) that can be found in Section~3.

Still, I had the strong feeling that this identity must be known.
So, if I cannot find it in classical sources, what about
``non-classical" sources?
I remembered that Alain Lascoux had
devoted one chapter of his book~\cite{\LascAZ} on symmetric
functions to orthogonal polynomials, revealing there that orthogonal
polynomials can be seen as Schur functions of quadratic shapes, and
demonstrating that formal identities for orthogonal polynomials 
can be conveniently established
by adopting this point of view. So I 
consulted \cite{\LascAZ}, and I quickly realised that
Proposition~8.4.1 in the book addresses~(\AA); it needed some more 
work 
to see what exactly was contained in that proposition.\footnote{Lascoux's
book is written in the (for many: foreign) language of plethystic operators
on symmetric polynomials. I therefore provide a translation of
the parts of \cite{\LascAZ} that are relevant to our discussion
into ``standard language" further below in Subsection~2.1.}

Lascoux attributes his proposition to
Christoffel, without any specific reference. 
With this information in hand, I returned to Szeg\H o's 
book \cite{\SzegAA} and made a text search for ``Christoffel".
I finally found the relevant theorem:
Theorem~2.5. I believe that, having a look at that theorem,
the reader will excuse myself that I
did not recognise that theorem as the one that is relevant for
our Theorem~{\TA} on my first attempt.\footnote{The reader may judge
her/himself: I present the theorem further below in Subsection~2.2, 
together with
explanations how this connects to our discussion.} 
In particular, \cite{\SzegAA, Theorem~2.5} does not
say anything about the proportionality factor between the two
sides in~(\AA) (as opposed to Lascoux, even if the expression he
gives is not correct; he does provide an argument though\footnote{Lascoux 
refers to ``the Bazin formula" without any reference. He may be
excused for that: a glance at the index of~\cite{\LascAZ} leads
one to \cite{\LascAZ, Lemma~A.1.1}.}).
Szeg\H o tells that \cite{\SzegAA, Theorem~2.5} is due to
Christoffel~\cite{\ChriAA}, but only in the special case of Legendre
polynomials
(indeed, at the end of \cite{\ChriAA} there appears Theorem~{\TA}
in that special case), a fact that also seems to have escaped many
researchers in the theory of orthogonal polynomials.

Eventually, I found that Theorem~\TA\ appears, correctly stated,
as Theorem~1 in the relatively recent article \cite{\ElouAA}
by Elouafi. However, the proof given there is incomplete.\footnote{The proof of
Lemma~4 in \cite{\ElouAA} only works for pairwise
distinct~$\al_i$'s. It is probably possible to complete the argument
even with an accordingly 
weakened version of that lemma. In our completion of Elouafi's
proof in Section~5 we prefer to complete the proof of the lemma.} 
I present a completion of this proof
in Section~5.

With all this knowledge, I consulted Mourad Ismail and asked
him if he knows the formula, respectively can refer me to a source
in the literature. He immediately pointed out that the right-hand
side determinant of~(\AA) features in ``Christoffel's theorem"
about the orthogonal polynomials with respect to the 
measure defined by the density $\prod _{i=1} ^{d-1}(x+x_i)\,d\mu(x)$ 
(with $d\mu(x)$ the density of the original orthogonality measure), 
that is, in \cite{\SzegAA,
Theorem~2.5} respectively \cite{\IsmaAA, Theorem~2.7.1}.
However, the conclusion of a longer discussion was that he had
not seen this formula earlier.

Finally, when I posted \cite{\KratCP} (containing
the extension of (\AA) to a rational deformation of the
density $d\mu(x)$) on the {\tt ar$\chi$iv}, Arno Kuijlaars
brought the article \cite{\BaDSAA} to my mind. Indeed, 
Equation~(2.6) in \cite{\BaDSAA} is equivalent to~(\AA), and it is
pointed out there that this result had been earlier obtained by
Br\'ezin and Hikami in \cite{\BrHiAA, Eq.~(14)}.
It requires some translational work to see this though, see Subsection~2.3. 

\medskip
In the next subsection, I provide a translation, into ``standard
English", of Lascoux's
rendering of Theorem~\TA. Then, in Subsection~2.2, I present
``Christoffel's theorem" and explain its connection to Theorem~\TA.
Finally, in Subsection~2.3 I translate the random matrix result
\cite{\BrHiAA, Eq.~(14)} into the language that we use here to see
that it is indeed equivalent to~(\AA).

\medskip
{\smc 2.1. Lascoux's Proposition~{\rm 8.4.1} in \cite{\LascAZ}.}
This proposition says
that, {\it given alphabets $\Bbb A$ and\/ $\Bbb B=\{b_1,b_2,\dots,b_{k+1}\}$,%
\footnote{In the statement of \cite{\LascAZ, Prop.~8.4.1},
the resultant $R(x,\Bbb B)$ must be replaced by the Vandermonde
product $\Delta(x+\Bbb B)$, as is done
in the proof of that proposition in \cite{\LascAZ}. 
I have done this correction here. 
Furthermore I simplified the statement
by incorporating the variable~$x$ in the alphabet~$\Bbb B$, which
means to replace ``$-\Bbb B-x$'' by ``$-\Bbb B$" and
``$x+\Bbb B$'' by ``$\Bbb B$". It is obvious that
Lascoux was aware of this simplification. However, he needed to
formulate the statement in that particular way in order to relate
it to the classical result \cite{\SzegAA, Theorem~2.5} (see also
\cite{\IsmaAA, Theorem~2.7.1}) 
in the theory of orthogonal polynomials.} 
$$
S_{(n+k)^n}(\Bbb A-\Bbb B)\Delta(\Bbb B)
$$
is proportional, up to a factor independent of\/ $\Bbb B$, to}
$$
\det_{1\le i,j\le k+1}\big(P_{n-1+j}(b_j)\big).
$$
The proportionality factor is given in the proof
of \cite{\LascAZ, Prop.~8.4.1} (except for the overall sign%
\footnote{Lascoux's attitude towards signs is best described by
himself: {\it``\dots \ with signs that specialists will know
how to write."} \cite{\LascAZ, comment added below Eq.~(3.1.5)}}),
but it has not been correctly worked out.

In order to understand the connection with Theorem~{\TA}, let me
translate Lascoux's language to the notation that I use here in
this paper. First of all, as mentioned in 
%??Musz gegebenfalls angepaszt werden.
Footnote~7, the symbol
$\Delta(\Bbb B)$
denotes the Vandermonde product (see \cite{\LascAZ, bottom of p.~11})
$$
\Delta(\Bbb B)=
\prod _{1\le i<j\le k+1} ^{}(b_j-b_i).
$$
Next, $S_{r^s}(\Bbb C)$ is the Schur function of rectangular shape
$r^s=(r,r,\dots,r)$ (with $s$~occurrences of~$r$) 
in the alphabet~$\Bbb C$ (not to be confused with 
the complex numbers!), defined by
(see \cite{\LascAZ, Eq.~(1.4.3)})
$$
S_{r^s}(\Bbb C)=\det_{1\le i,j\le s}\big(S^{r+j-i}(\Bbb C)\big)
=(-1)^{\binom s2}\det_{1\le i,j\le s}\big(S^{r+i+j-s-1}(\Bbb C)\big),
$$
where $S^a(\Bbb C)$ is the complete homogeneous symmetric function
of degree~$a$ in the alphabet~$\Bbb C$. Thus,
$$
\align
S_{(n+k)^n}(\Bbb A-\Bbb B)&=(-1)^{\binom n2}
\det_{1\le i,j\le n}\big(S^{(n+k)+i+j-n-1}(\Bbb A-\Bbb B)\big)\\
&=(-1)^{\binom n2}\det_{0\le i,j\le n-1}\big(S^{i+j+k+1}(\Bbb A-\Bbb B)\big).
\endalign
$$
Here (this is implicit on page~6 of \cite{\LascAZ}),
$$
S^{m}(\Bbb A-\Bbb B)=\sum_{a=0}^{m}S^{m-a}(\Bbb A)S^{a}(-\Bbb B)
=\sum_{a=0}^{m}(-1)^aS^{m-a}(\Bbb A)\La^{a}(\Bbb B),
$$
where $\La^{a}(\Bbb B)$ is the elementary symmetric function of
degree~$a$ in $b_1,b_2,\dots,b_{k+1}$. Furthermore (cf\. \cite{\LascAZ,
second display on p.~115}), 
we have $S^m(\Bbb A)=\mu_m$, with $\mu_m$ being ``our" $m$-th
moment of the linear functional~$L$. Thus, we recognise that
$S_{(n+k)^n}(\Bbb A-\Bbb B)$ is, up to the sign $(-1)^{\binom n2}$, 
in our notation the Hankel determinant 
$$
\det\limits_{0\le i,j\le n-1}\left(
\mu^{i+j}\prod _{\ell=1} ^{k+1}(-b_\ell+\mu)\right)
$$
on the left-hand side of~(\AA) (with $x_\ell$ replaced by $-b_\ell$
and $d$ replaced by $k+1$).

On the other hand, Lascoux's polynomials $P_n(x)$ are {\it orthonormal\/} with
respect to the linear functional~$L$ with moments $\mu_m$, $m=0,1,\dots$,
and are given by (cf\. \cite{\LascAZ, Theorem~8.1.1})
$$
P_n(x)=\frac {1} {\big((-1)^nS_{(n-1)^n}(\Bbb A)S_{n^{n+1}}(\Bbb A)\big)^{1/2}} 
S_{n^n}(\Bbb A-x), \quad \text{for }n\ge0,
$$
while ``our" {\it monic} orthogonal polynomials $p_n(x)$ are given by
(cf\. \cite{\LascAZ, Theorem 8.1.1 and second display on
p.~116})
$$
p_n(x)=\frac {1} {(-1)^nS_{(n-1)^n}(\Bbb A)} 
S_{n^n}(\Bbb A-x), \quad \text{for }n\ge0.
$$
Thus, we see that Lascoux's determinant $\det_{1\le i,j\le
k}\big(P_{n-1+j}(b_j)\big)$ is, up to some overall factor,
``our" determinant $\det_{1\le i,j\le
k}\big(p_{n-1+j}(b_j)\big)$
on the right-hand side of~(\AA) (with $x_\ell$ replaced by $-b_\ell$
and $d$ replaced by $k+1$).

It should now be clear to the reader 
that Lascoux's Proposition~8.4.1 in \cite{\LascAZ} is equivalent to
Theorem~{\TA}, except that he did not bother to figure out the
correct sign, and that he did not get the proportionality factor
right (both of which being very understandable
given the complexity of the task \dots; in fact,
in order to not risk to also fail, I do not attempt to present
the correct proportionality factor or sign in Lascoux's notation ---
in ``standard" notation, the correct identity is given in~(\AA).)

\medskip
{\smc 2.2. ``Christoffel's Theorem".} This theorem
(cf\. \cite{\SzegAA, Theorem~2.5} or
\cite{\IsmaAA, Theorem~2.7.1}) says: {\it in the setting of Section~{\rm1},
the polynomials (in~$x$)
$$
\frac {\displaystyle\det_{1\le i,j\le d}\pmatrix 
p_{n+i-1}(x_j),&1\le j\le d-1\\
p_{n+i-1}(x),&j=d\hfill\endpmatrix}
{
\prod\limits _{i=1} ^{d-1}(x-x_i)}, \quad n=0,1,\dots,
\tag\AZa$$
form a sequence of orthogonal polynomials with respect to the
linear functional with moments
$$
\mu^{n}\prod _{\ell=1} ^{d-1}(\mu-x_\ell), \quad n=0,1,\dots.
\tag\AZb$$}%

The reader may now understand why I did not notice that Theorem~\TA\
is hidden in the above assertion on my first attempt to find
it in Szeg\H o's book~\cite{\SzegAA}. 
I did of course see that the determinant in~(\AZa)
is our determinant on the right-hand side of~(\AA) (with $x_d=x$),
and that the moments in~(\AZb) are ``almost" the entries in the Hankel
determinant on the left-hand side of~(\AA). There are however two obstacles
to overcome in order to ``extract"~(\AA) out of the above assertion:
first, one has to recall a certain determinant formula for the
orthogonal polynomials with respect to a given moment sequence, namely
Lemma~\TD\ in Section~4. Applied to the moments in~(\AZb), it
produces indeed the Hankel determinant on the left-hand side
of~(\AA). The uniqueness of orthogonal polynomials up to scaling
then implies that the determinants on the left-hand side of~(\AA) 
(with $x_d=x$)
and the expression~(\AZa)
agree up to a multiplicative constant
(meaning: independent of $x=x_d$). So, second, this constant
has to be computed. A researcher in the theory of orthogonal
polynomials does not really care about the normalisation of the
orthogonal polynomials and this seems to be the reason why 
apparently nobody has done it there, although this is not really
difficult, see Sections~4 and~5 for two slightly different arguments.

\medskip
{\smc 2.3. Expectation of a product of characteristic polynomials
of random Hermitian matrices.} 
Let $d\mu(u)$ be the density of some positive measure with infinite
support all of whose moments exist.
Equation~(14) in \cite{\BrHiAA} (cf\. also \cite{\BaDSAA, Eq.~(2.6)}) reads
$$
\left\langle
\prod _{j=1} ^{d}D_n(x_j,H)
\right\rangle_\mu
=\frac {\det_{1\le i,j\le d}\big(p_{n+i-1}(x_j)\big)} {
\prod _{1\le i<j\le d} ^{}(x_j-x_i)}.
\tag\AZc$$
(I have changed notation so that it is in line with our notation.)
Here, the left-hand side is an expectation for products of characteristic
polynomials of random Hermitian matrices. However, this fact
does not need to concern us. According to \cite{\BrHiAA, Eq.~(5)}
(see also \cite{\BaDSAA, Eq.~(1.3)})
it can be expressed as
$$\multline
\left\langle
\prod _{j=1} ^{d}D_n(x_j,H)
\right\rangle_\mu
\\
=
\frac {1} {Z_n}\int\cdots\int
\bigg(\prod _{j=1} ^{d}
\prod _{i=1} ^{n}(x_j-u_i)
\bigg)
\bigg(
\prod _{1\le i<j\le n} (u_j-u_i)
\bigg)^2\,
d\mu(u_1)\,d\mu(u_2)\cdots d\mu(u_n),
\endmultline$$
where
$$
Z_n=\int\cdots\int
\bigg(
\prod _{1\le i<j\le n} (u_j-u_i)
\bigg)^2\,
d\mu(u_1)\,d\mu(u_2)\cdots d\mu(u_n).
$$ 
Now, by Heine's formula (cf\. \cite{\SzegAA, Eq.~(2.2.11)},
\cite{\IsmaAA, Cor.~2.1.3}, or Lemma~\TFa), we have
$$
Z_n=n!\,\det_{0\le i,j\le n-1}(\mu_{i+j})
$$
and
$$
\left\langle
\prod _{j=1} ^{d}D_n(x_j,H)
\right\rangle_\mu
=n\,!
\det_{0\le i,j\le n-1}\left(\int 
u^{i+j}\bigg(\prod _{j=1} ^{d}
(x_j-u)\bigg)\,d\mu(u)\right).
$$
Since
$$
\int u^{i+j}
\bigg(\prod _{j=1} ^{d}
(x_j-u)\bigg)\,d\mu(u)=
\mu^{i+j}\prod _{j=1} ^{d}
(x_j-\mu)
$$
(again using umbral notation), the equivalence of (\AZc) and (\AA) now becomes
obvious.

\subhead 3. First proof of Theorem \TA\ --- condensation\endsubhead
In this section, we present the author's proof of Theorem~\TA, which uses
the method of condensation (frequently referred to as ``Dodgson
condensation"). This method provides inductive proofs
that are based on a determinant identity due to Jacobi (see
Proposition~\TB\ below).

\medskip
For convenience, we change notation slightly. Instead of the
polynomials $p_n(x)$, let us consider
the polynomials $f_n(x)$ defined by 
$$
f_n(x)=(x+s_{n-1})f_{n-1}(x)-t_{n-2}f_{n-2}(x),
\quad \text{for }n\ge1,
\tag\BA
$$
with $f_0(x) = 1$ and $f_{-1}(x)=0$. (It should be noted that the only
difference between the recurrences (\AB) and (\BA) is the sign in
front of~$s_{n-1}$.) Using these polynomials, the formula~(\AA) 
can be rewritten as
$$
\frac {\det\limits_{0\le i,j\le n-1}\left(
\mu^{i+j}\prod _{\ell=1} ^{d}(x_\ell+\mu)\right)}
{\det\limits_{0\le i,j\le n-1}\left(\mu_{i+j}\right)}
=
\frac {\displaystyle\det_{1\le i,j\le d}\left(f_{n+i-1}(x_j)\right)}
{
\prod\limits _{1\le i<j\le d} ^{}(x_j-x_i)}.
\tag\BB
$$
(That is, we got rid of the signs on the right-hand side of~(\AA).)
Our proof of (\BB) will be based on the {\it method of condensation}
(see \cite{\KratBN, Sec.~2.3}). The ``backbone" of this method is
the following determinant identity due to Jacobi.

\proclaim{Proposition \TB}
Let $A$ be an $N\times N$ matrix. Denote the submatrix of $A$ in which
rows $i_1,i_2,\dots,i_k$ and columns $j_1,j_2,\dots,j_k$ are 
omitted by $A_{i_1,i_2,\dots,i_k}^{j_1,j_2,\dots,j_k}$. Then we have
$$
\det A\cdot \det A_{i_1,i_2}^{j_1,j_2}=\det A_{i_1}^{j_1}\cdot 
\det A_{i_2}^{j_2}-
\det A_{i_1}^{j_2}\cdot \det A_{i_2}^{j_1}
\tag\BC
$$
for all integers $i_1,i_2,j_1,j_2$ with $1\le i_1<i_2\le N$
and $1\le j_1<j_2\le N$.
\endproclaim

The second ingredient of the proof of (\BB) will be the 
Hankel determinant identity below, which, as its proof will reveal, 
is actually a consequence of the condensation formula in~(\BC).

\proclaim{Lemma \TC}
Let $(c_n)_{n\ge0}$ be a given sequence, and $\al$ and $\be$ be
variables. Then, for all positive integers $n$, we have
$$\multline 
(\be-\al)\det_{0\le i,j\le n-1}\big(\al\be
c_{i+j}+(\al+\be)c_{i+j+1}+c_{i+j+2}\big)
\det_{0\le i,j\le n}\big(c_{i+j}\big)\\
=
\det_{0\le i,j\le n-1}\big(\al c_{i+j}+c_{i+j+1}\big)
\det_{0\le i,j\le n}\big(\be c_{i+j}+c_{i+j+1}\big)\\
-
\det_{0\le i,j\le n-1}\big(\be c_{i+j}+c_{i+j+1}\big)
\det_{0\le i,j\le n}\big(\al c_{i+j}+c_{i+j+1}\big).
\endmultline
\tag\BD$$
\endproclaim

\demo{Proof}
By using multilinearity in the rows, it is easy to see
(cf\. also \cite{\KratCI, Lemma~4}) that
$$\det_{0\le i,j\le M}(\al c_{i+j}+c_{i+j+1})=
\sum _{r=0} ^{M+1}\al^r
\det_{0\le i,j\le M}(c_{i+j+\chi(i\ge r)}),
\tag\BE$$
where $\chi(\Cal S)=1$ if $\Cal S$ is
true and $\chi(\Cal S)=0$ otherwise.
If we apply this identity to the first determinant
on the left-hand side of (\BD), then we obtain
$$\align 
\det_{0\le i,j\le n-1}\big(\al\be
c_{i+j}+(\al+\be)&c_{i+j+1}+c_{i+j+2}\big)\\
&=
\det_{0\le i,j\le n-1}\big(\al(\be
c_{i+j}+c_{i+j+1})+(\be c_{i+j+1}+c_{i+j+2})\big)\\
&=
\sum_{r=0}^n\al^r\det_{0\le i,j\le n-1}\big(
\be c_{i+j+\chi(i\ge r)}+c_{i+j+\chi(i\ge r)+1}\big).
\endalign$$
Now we use multilinearity in rows $0,1,\dots,r-1$ and in rows
$r,r+1,\dots,n-1$ separately. This leads to
$$\align 
\det_{0\le i,j\le n-1}\big(\al\be
c_{i+j}+&(\al+\be)c_{i+j+1}+c_{i+j+2}\big)\\
&=
\sum_{r=0}^n\al^r
\sum_{s=0}^{r}\be^s
\sum_{t=r}^n\be^{t-r}
\det_{0\le i,j\le n-1}\big(
c_{i+j+\chi(i\ge s)+\chi(i\ge t)}\big)\\
&=
\sum_{0\le s\le t\le n}
\bigg(\sum_{r=s}^t \al^r\be^{s+t-r}\bigg)
\det_{0\le i,j\le n-1}\big(
c_{i+j+\chi(i\ge s)+\chi(i\ge t)}\big)\\
&=
\frac 1{\be-\al}
\sum_{0\le s\le t\le n}
\left(\al^s\be^{t+1}-\al^{t+1}\be^{s}\right)
\det_{0\le i,j\le n-1}\big(
c_{i+j+\chi(i\ge s)+\chi(i\ge t)}\big).
\endalign$$
We substitute this as well as (\BE) (wherever it can be applied)
in~(\BD). The factor $\be-\al$ cancels. Subsequently, we
compare coefficients of $\al^s\be^{t+1}$, respectively of
$\al^{t+1}\be^s$, for $0\le s\le t\le n$. Thus we see that we need to show
$$\multline 
\det_{0\le i,j\le n-1}\big(
c_{i+j+\chi(i\ge s)+\chi(i\ge t)}\big)
\det_{0\le i,j\le n}\big(
c_{i+j}\big)\\
=
\det_{0\le i,j\le n-1}\big(c_{i+j+\chi(i\ge s)}\big)
\det_{0\le i,j\le n}\big(c_{i+j+\chi(i\ge t+1)}\big)\\
-
\det_{0\le i,j\le n-1}\big(c_{i+j+\chi(i\ge t+1)}\big)
\det_{0\le i,j\le n}\big(c_{i+j+\chi(i\ge s)}\big),
\endmultline
\tag\BF$$
and this would moreover be sufficient for the proof of~(\BD).
As it turns out, the choice of
$$
A=\pmatrix 
c_0&c_1&c_2& \dots&c_n&0\\
c_1&c_2&c_3& \dots&c_{n+1}&0\\
c_2&c_3&c_4& \dots&c_{n+2}&0\\
\hdotsfor6\\
c_n&c_{n+1}&c_{n+2}& \dots&c_{2n}&0\\
c_{n+1}&c_{n+2}&c_{n+3}& \dots&c_{2n+1}&1
\endpmatrix
$$
and $i_1=s$, $i_2=t+1$, $j_1=n$, $j_2=n+1$ in Proposition~{\TB} yields
exactly~(\BF).\quad \quad \qed
\enddemo

We are now ready for the proof of Theorem~\TA, which, as we have
seen, is equivalent to~(\BB).

\demo{Proof of (\BB)}
We prove (\BB), in the form
$$
 {\displaystyle\det_{1\le i,j\le d}\left(f_{n+i-1}(x_j)\right)}
=\bigg(
\prod\limits _{1\le i<j\le d} ^{}(x_j-x_i)\bigg)
\frac {\det\limits_{0\le i,j\le n-1}\left(
\mu^{i+j}\prod _{\ell=1} ^{d}(x_\ell+\mu)\right)}
{\det\limits_{0\le i,j\le n-1}\left(\mu_{i+j}\right)},
\tag\BG
$$
by induction on $d$. For $d=0$, the formula is trivially true.
For $d=1$, the formula has been proven by 
Mu, Wang and Yeh in \cite{\MuWYAA, Theorem~1.3}
in a different but equivalent form (see also \cite{\CJKrAB,
Eq.~(3.2)}).

Let $\LHS_{d,n}(x_1,\dots,x_d)$ denote the left-hand side of (\BG).
For the induction step, we observe that, according to~(\BC) with $N=d$,
$A=\big(f_{n+i-1}(x_j)\big)_{1\le i,j\le d}$, $i_1=j_1=1$ and
$i_d=j_d=d$, we have
$$\multline 
\LHS_{d,n}(x_1,\dots,x_d)\LHS_{d-2,n+1}(x_2,\dots,x_{d-1})\\
=
\LHS_{d-1,n}(x_1,\dots,x_{d-1})\LHS_{d-1,n+1}(x_2,\dots,x_d)\\
-\LHS_{d-1,n}(x_2,\dots,x_d)\LHS_{d-1,n+1}(x_1,\dots,x_{d-1}).
\endmultline
\tag\BH$$
This can be seen as a recurrence formula for
$\LHS_{d,n}(x_1,\dots,x_d)$,
as one can use it to express $\LHS_{d,n}(x_1,\dots,x_d)$ in terms
of expressions of the form $\LHS_{e,m}(x_a,\dots,x_b)$ with~$e$
smaller than~$d$. Hence, for the proof of (\BG) it suffices to
verify that the right-hand side of~(\BG) satisfies the same
recurrence. Consequently, we substitute this right-hand side in~(\BH).
After cancellation of factors that are common to both sides, we
arrive at
$$\multline 
(x_d-x_1)
\det\limits_{0\le i,j\le n-1}\left(
\mu^{i+j}\prod _{\ell=1} ^{d}(x_\ell+\mu)\right)
\det\limits_{0\le i,j\le n}\left(
\mu^{i+j}\prod _{\ell=2} ^{d-1}(x_\ell+\mu)\right)\\
=
\det\limits_{0\le i,j\le n-1}\left(
\mu^{i+j}\prod _{\ell=1} ^{d-1}(x_\ell+\mu)\right)
\det\limits_{0\le i,j\le n}\left(
\mu^{i+j}\prod _{\ell=2} ^{d}(x_\ell+\mu)\right)\\
-
\det\limits_{0\le i,j\le n-1}\left(
\mu^{i+j}\prod _{\ell=2} ^{d}(x_\ell+\mu)\right)
\det\limits_{0\le i,j\le n}\left(
\mu^{i+j}\prod _{\ell=1} ^{d-1}(x_\ell+\mu)\right).
\endmultline$$
This is the special case of Lemma~{\TC} where $c_i=\mu^{i}\prod
_{\ell=2} ^{d-1}(x_\ell+\mu)$, $\al=x_1$ and $\be=x_d$.\quad \quad \qed
\enddemo

\subhead 4. Second proof of Theorem \TA\ --- theory of orthogonal polynomials\endsubhead
Here we describe a proof of Theorem~\TA\ that is based on facts from
the theory of orthogonal polynomials. We follow largely Lascoux's
arguments in the proof of Proposition~8.4.1 in~\cite{\LascAZ}.
They show that, using the uniqueness up to scaling of orthogonal polynomials
with respect to a given linear functional,
the right-hand side and the left-hand side in~(\AA) agree up
to a multiplicative constant. For the determination of this
constant we provide a simpler argument than the one given in~\cite{\LascAZ}.

\medskip
We prove (\AA) in the form
$$
\frac {\det\limits_{0\le i,j\le n-1}\left(
\mu^{i+j}\prod _{\ell=1} ^{d}(\mu-x_\ell)\right)}
{\det\limits_{0\le i,j\le n-1}\left(\mu_{i+j}\right)}
=(-1)^{nd}
\frac {\displaystyle\det_{1\le i,j\le d}\left(p_{n+i-1}(x_j)\right)}
{
\prod\limits _{1\le i<j\le d} ^{}(x_j-x_i)}.
\tag\CA
$$
It should be recalled that, with $L$ denoting the functional of
orthogonality for the polynomials $\big(p_n(x))_{n\ge0}$, we have
$L(x^n)=\mu_n$, where we still use the umbral notation $\mu^n\equiv\mu_n$.

We start with a classical fact from the theory of orthogonal polynomials
(cf\. \cite{\SzegAA, Eq.~(2.2.9)} or \cite{\IsmaAA, Eq.~(2.1.10)}).

\proclaim{Lemma \TD}
Let $M$ be a linear functional on polynomials in $x$ with
moments~$\nu_n$, $n=0,1,\dots$, such that all Hankel determinants 
$\det_{0\le i,j\le n}(\nu_{i+j})$, $n=0,1,\dots$, 
are non-zero. 
Then the determinants
$$
\det_{0\le i,j\le n-1}\left(\nu_{i+j+1}-\nu_{i+j}x\right)
$$
are a sequence of orthogonal polynomials with respect to~$M$.
\endproclaim

\demo{Proof}
We have
$$\align
\det_{0\le i,j\le n-1}&\left(\nu_{i+j+1}-\nu_{i+j}x\right)\\
&=\det\pmatrix 
1&0&0&\dots&0\\
\nu_0&\nu_1-\nu_0x&\nu_2-\nu_1x&\dots&\nu_n-\nu_{n-1}x\\
\nu_1&\nu_2-\nu_1x&\nu_3-\nu_2x&\dots&\nu_{n+1}-\nu_{n}x\\
\hdotsfor5\\
\nu_{n-1}&\nu_n-\nu_{n-1}x&\nu_{n+1}-\nu_nx&\dots&\nu_{2n-1}-\nu_{2n-2}x
\endpmatrix\\
&=
\det\pmatrix 
1&x&x^2&\dots&x^n\\
\nu_0&\nu_1&\nu_2&\dots&\nu_n\\
\nu_1&\nu_2&\nu_3&\dots&\nu_{n+1}\\
\hdotsfor5\\
\nu_{n-1}&\nu_n&\nu_{n+1}&\dots&\nu_{2n-1}
\endpmatrix.
\tag\CAa
\endalign$$
It is straightforward to check that the determinant in the last line
is orthogonal with respect to $1,x,x^2,\dots,x^{n-1}$.
Moreover, the coefficient of~$x^n$ is 
$\pm\det_{0\le i,j\le n-1}(\nu_{i+j})$, which by
assumption is non-zero so that the determinant in the assertion of the
lemma is a polynomial of
degree~$n$, as desired.\quad \quad \qed
\enddemo

\remark{Remark}
The determinant in the last line of (\CAa) represents another classical
determinantal formula for orthogonal polynomials expressed in terms of
the moments of the corresponding linear functional of orthogonality,
see \cite{\SzegAA, Eq.~(2.2.6)} or \cite{\IsmaAA, Eq.~(2.1.11)}.
\endremark

\demo{Proof of (\CA)}
Using Lemma~\TD\ with $\nu_n=\mu^{n}\prod _{\ell=1} ^{d-1}(\mu-x_\ell)$,
we see that the determinants in the numerator of the left-hand side of~(\CA),
$$
\det\limits_{0\le i,j\le n-1}\left(
\mu^{i+j}\prod _{\ell=1} ^{d}(\mu-x_\ell)\right),
$$
seen as polynomials in $x_d$, are a sequence of orthogonal
polynomials for the linear functional with moments
$$
\mu^{n}\prod _{\ell=1} ^{d-1}(\mu-x_\ell),\quad 
n=0,1,\dots.
\tag\CB
$$
Clearly,
in terms of the functional $L$ (now acting on polynomials in~$x_d$)
of orthogonality for the polynomials
$\big(p_n(x_d)\big)_{n\ge0}$,
this linear functional with moments~(\CB) can be expressed as
$$
p(x_d)\mapsto L\bigg(p(x_d)\cdot 
\prod _{\ell=1} ^{d-1}(x_d-x_\ell)\bigg).
\tag\CC
$$

We claim that also the right-hand side of (\CA) gives a sequence
of orthogonal polynomials (in~$x_d$) with respect to the linear
functional~(\CC). The first (and easy) observation is that
the right-hand side of~(\CA) has indeed degree~$n$ in~$x_d$.

Let us denote the right-hand side
of (\CA) by $q_n(x_d)$. When we apply the functional~(\CC) to 
$x_d^sq_n(x_d)$, for $0\le s\le n-1$,
then, up to factors which are independent of~$x_d$, we obtain
$$
L\left(x_d^s\det_{1\le i,j\le d}\left(p_{n+i-1}(x_j)\right)\right).
$$
By expanding the determinant with respect to the last column, this
becomes a linear combination of terms of the form $L(x_d^sp_{n+i-1}(x_d))$. 
Since $i\ge1$ and $s\le n-1$, all of them vanish, proving our claim.

By symmetry, the same argument can also be made for any $x_\ell$
with $1\le \ell\le d-1$.

The fact that orthogonal polynomials with respect to a particular
linear functional are unique up to multiplicative constants then
implies that 
$$
 {\det\limits_{0\le i,j\le n-1}\left(
\mu^{i+j}\prod _{\ell=1} ^{d}(\mu-x_\ell)\right)}
=C
\frac {\displaystyle\det_{1\le i,j\le d}\left(p_{n+i-1}(x_j)\right)}
{
\prod\limits _{1\le i<j\le d} ^{}(x_j-x_i)},
\tag\CD$$
where $C$ is independent of the variables $x_1,x_2,\dots,x_d$.
In order to compute~$C$, we divide both sides by
$x_1^nx_2^n\cdots x_d^n$, and then compute the limits as
$x_d\to\infty$, \dots
$x_2\to\infty$,
$x_1\to\infty$, in this order. It is not difficult to see that in
this manner the above equation reduces to
$$
 {\det\limits_{0\le i,j\le n-1}\left(
\mu_{i+j}(-1)^d\right)}
=C
\det A,
$$
where $A$ is a lower triangular matrix with ones on the diagonal.
Hence, we get $C=(-1)^{nd}\det_{0\le i,j\le n-1}(\mu_{i+j})$,
as desired.\quad \quad \qed
\enddemo

\subhead 5. Third proof of Theorem \TA\ --- vanishing of polynomials\endsubhead
The purpose of this section is to present a completed version of
Elouafi's proof of Theorem~\TA\ in~\cite{\ElouAA}. It is based
on a vanishing argument: it is shown that the left-hand side
of~(\AA) vanishes if and only if the right-hand side does.
Since both sides have the same leading monomial as polynomials
in the $x_i$'s, it follows that they must be the same up to
a multiplicative constant. This constant is then determined in
the last step.

\medskip
To begin with, we need some preparations.
Let us write
$$
R(x_1,x_2,\dots,x_d):=\frac {\displaystyle\det_{1\le i,j\le d}\left(p_{n+i-1}(-x_j)\right)}
{
\prod\limits _{1\le i<j\le d} ^{}(x_i-x_j)}
$$
for the expression on the right-hand side of (\AA), forgetting the
sign.
Since the numerator is skew-symmetric in the $x_i$'s, it is
divisible by the Vandermonde product $\prod _{1\le i<j\le d}
^{}(x_i-x_j)$ in the denominator, so that $R(x_1,x_2,\dots,x_d)$ 
is actually a (symmetric) polynomial in the $x_i$'s. Thus, while in its
definition it seems problematic to substitute the same value
for two different $x_i$'s in $R(x_1,x_2,\dots,x_d)$, this is
actually not a problem. Nevertheless, it would also be good to have
an explicit form for such a case available as well. This is afforded
by the proposition below.

\proclaim{Proposition \TEa}
Let $y_1,y_2,\dots,y_e$ be variables and $m_1,m_2,\dots,m_e$ be
non-negative integers with $m_1+m_2+\dots+m_e=d$. Then we have
$$
R(y_1,\dots,y_1,y_2,\dots,y_2,\dots,y_e,\dots,y_e)
=\frac {\det M_{m_1,m_2,\dots,m_e}(y_1,y_2,\dots,y_e)} 
{\prod\limits _{1\le i<j\le e} ^{}(y_i-y_j)^{m_im_j}},
$$
where $y_i$ is repeated $m_i$ times in the argument of $R$ on the
left-hand side. The matrix $M_{m_1,m_2,\dots,m_e}(y_1,y_2,\dots,y_e)$
is defined by
$$
M_{m_1,m_2,\dots,m_e}(y_1,y_2,\dots,y_e)
=\big(M_{m_1}(y_1)\ M_{m_2}(y_2)\ \cdots\ M_{m_e}(y_e)\big),
$$
where $M_m(y)$ is the $d\times m$ matrix
$$
M_m(y)=\pmatrix \dfrac {p^{(j-1)}_{n+i-1}(-y)}
{(j-1)!}\endpmatrix_{1\le i\le d,\ 1\le j\le m},
$$
with $p^{(j)}_n(y)$ denoting the $j$-th derivative of $p_n(y)$ with respect
to~$y$.
\endproclaim

\demo{Proof}
We have to compute the limit
$$
\lim_{x_1\to y_1,\dots,x_{m_1}\to y_1,\dots,x_{m_1+\dots+m_{e-1}+1}\to
y_e,\dots,x_d\to y_e}R(x_1,x_2,\dots,x_d).
$$
Since we know that $R(x_1,x_2,\dots,x_d)$ is in fact a polynomial in
the~$x_i$'s, we have a large flexibility of how to compute this limit.
We choose to do it as follows: 
we put $x_i=y_1+ih$ for $i=1,2,\dots,m_1$,
$x_i=y_2+ih$ for $i=m_1+1,m_1+2,\dots,m_1+m_2$,
etc. In the end we let $h$ tend to zero.

We describe how this works for the first group of variables,
for the other the procedure is completely analogous.
In the matrix appearing in the numerator of the definition of
$R(x_1,x_2,\dots,x_d)$, we replace column~$j$ by
$$
\sum_{k=1}^{j}(-1)^{j-k}\binom {j-1}{k-1}(\text{column }k),\quad 
\text{for }j=1,2,\dots,m_1.
$$
Clearly, this modification of the matrix can be achieved by elementary
column operations, so that the determinant is not changed. Thus, we
obtain
$$
R(x_1,x_2,\dots,x_d)
=\frac {\det N_1}
{
\prod\limits _{1\le i<j\le d} ^{}(x_i-x_j)},
\tag\DAd$$
where
$$
\multline
N_1=\bigg(p_{n+i-1}(-x_1)\quad
p_{n+i-1}(-x_2)-p_{n+i-1}(-x_1)\quad \dots\quad \\
\sum_{k=1}^{m_1}(-1)^{m_1-k}\binom {m_1-1}{k-1}p_{n+i-1}(-x_k)\quad 
 \dots
\bigg)_{1\le i\le d}.
\endmultline
$$
(Here, the terms describe the columns of the matrix~$N_1$.)
We now perform the earlier described assignments for
$x_1,x_2,\dots,x_{m_1}$. Under these assignments, we have
$$
\prod\limits _{1\le i<j\le m_1} ^{}(x_i-x_j)
=
\prod\limits _{1\le i<j\le m_1} ^{}h(i-j)
=(-1)^{\binom {m_1}2}h^{\binom {m_1}2}
\prod\limits _{l=1}^{m_1-1}l!
$$
and
$$
\lim_{h\to0}\frac {1} {h^{j-1}}
\sum_{k=1}^{j}(-1)^{j-k}\binom {j-1}{k-1}p_{n+i-1}(-y_1-kh)
=(-1)^{j-1}p^{(j-1)}_{n+i-1}(-y_1).
$$
Therefore we get
$$
R(y_1,\dots,y_1,x_{m_1+1},\dots,x_d)
=\frac {\det N_2}
{
\bigg(\prod\limits _{l=1}^{m_1-1}l!\bigg)\bigg(
\prod\limits _{j=m_1+1} ^{d}(y_1-x_j)^{m_1}\bigg)\bigg(
\prod\limits _{m_1+1\le i<j\le d} ^{}(x_i-x_j)\bigg)},
$$
where
$$
N_2=\bigg(p_{n+i-1}(-y_1)\quad
p^{(1)}_{n+i-1}(-y_1)\quad \dots\quad 
p^{(m_1-1)}_{n+i-1}(-y_1)\quad 
p_{n+i-1}(-x_{m_1+1})\quad 
 \dots
\bigg)_{1\le i\le d}.
$$

To finish the argument, one has to proceed analogously for the
remaining groups of~$x_i$'s and finally put the arising factorials
into the columns of the determinant.\quad \quad \qed
\enddemo

The following auxiliary result is \cite{\ElouAA, Lemma~3}.

\proclaim{Lemma \TE}
Let $q(x)=
\prod _{i=1} ^{d}(x+x_i)$. Furthermore, as before, we write $L$ for
the linear functional with respect to which the sequence
$\big(p_n(x)\big)_{n\ge0}$ is orthogonal. Then
$$
\det\limits_{0\le i,j\le n-1}\left(
\mu^{i+j}\prod _{\ell=1} ^{d}(x_\ell+\mu)\right)
=
\det\limits_{0\le i,j\le n-1}\left(
L\big(q(x)p_i(x)p_j(x)\big)\right).
\tag\DAa$$
\endproclaim

\demo{Proof}
We first rewrite the determinant on the left-hand side,
$$
\det\limits_{0\le i,j\le n-1}\left(
\mu^{i+j}\prod _{\ell=1} ^{d}(x_\ell+\mu)\right)
=
\det\limits_{0\le i,j\le n-1}\left(
L\big(x^{i+j}q(x)\big)\right).
\tag\DAb
$$
Setting $p_m(x)=\sum_{a=0}^mc_{ma}x^a$, for some
coefficients~$c_{ma}$, we may rewrite the right-hand side of (\DAa) as
$$
\det\limits_{0\le i,j\le n-1}\left(
L\big(q(x)p_i(x)p_j(x)\big)\right)
=
\det\limits_{0\le i,j\le n-1}\left(
\sum_{a=0}^i
\sum_{b=0}^jc_{ia}c_{jb}
L\big(q(x)x^{a+b}\big)\right).
$$
Consequently, letting $c_{ma}=0$ whenever $m<a$, we have
$$
\det\limits_{0\le i,j\le n-1}\left(
L\big(q(x)p_i(x)p_j(x)\big)\right)
=
\det\left(A \cdot B\cdot C\right),
\tag\DAc
$$
where
$$
A=(c_{ia})_{0\le i,a\le n-1},\quad 
B=\big(L\big(q(x)x^{a+b}\big)\big)_{0\le a,b\le n-1},\quad 
C=(c_{jb})_{0\le b,j\le n-1}.
$$
Since $A$ and $C$ are triangular matrices with ones on the main
diagonal (the latter follows from our assumption that the polynomials
$p_m(x)$ are monic), we conclude that $\det A=\det C=1$.
The combination of~(\DAb) and~(\DAc) then establishes the desired
assertion.\quad \quad \qed
\enddemo

The following lemma is \cite{\ElouAA, Lemma~4}, which appears there
with an incomplete proof (cf\. 
%??Musz gegebenfalls angepaszt werden.
Footnote~6).

\proclaim{Lemma \TF}
As polynomial functions in the variables $x_1,x_2,\dots,x_d$, 
the expression
$$
R(x_1,x_2,\dots,x_d)=
\frac {\displaystyle\det_{1\le i,j\le d}\left(p_{n+i-1}(-x_j)\right)}
{
\prod\limits _{1\le i<j\le d} ^{}(x_i-x_j)}
\tag\DA$$
vanishes in an extension field $\widehat K$ of the ground field~$K$ 
if and only if the determinant
$$
S(x_1,x_2,\dots,x_d)=
\det\limits_{0\le i,j\le n-1}\left(
\mu^{i+j}\prod _{\ell=1} ^{d}(x_\ell+\mu)\right)
\tag\DB
$$
vanishes.
\endproclaim

\demo{Proof}
Let $R(x_1,x_2,\dots,x_d)=0$, for some choice
of elements $x_i$ in~$\widehat K$. We assume that
$$\{x_1,x_2,\dots,x_d\}=\{y_1,y_2,\dots,y_e\},$$ 
with the $y_i$'s being
pairwise distinct, and where $y_i$ appears with multiplicity $m_i$
among the $x_j$'s, $i=1,2,\dots,e$. Since both $R(x_1,x_2,\dots,x_d)$
and $S(x_1,x_2,\dots,x_d)$ are symmetric polynomials in the $x_j$'s,
without loss of generality we may assume that $x_1=\dots=x_{m_1}=y_1$,
$x_{m_1+1}=\dots=x_{m_1+m_2}=y_2$, \dots, 
$x_{m_1+\dots+m_{e-1}+1}=\dots=x_{d}=y_e$.
Thus, in view of Proposition~\TEa, the vanishing of
$$
R(y_1,\dots,y_1,y_2,\dots,y_2,\dots,y_e,\dots,y_e),
$$
where the $y_i$'s are pairwise distinct and $y_i$ is repeated $m_i$
times, is equivalent to
the rows of the matrix $M_{m_1,m_2,\dots,m_e}(y_1,y_2,\dots,y_e)$
being linearly dependent. In other words, there exist constants
$c_i\in\widehat K$,
$i=0,1,\dots,n-1$, not all of them zero, such that
$$
\sum_{i=1}^{d}c_ip^{(j-1)}_{n+i-1}(-y_k)=0,\quad 
\text{for }k=1,2,\dots,e\text{ and }j=1,2,\dots,m_k.
\tag\DC
$$
Define $g(x):=\sum_{i=1}^{d}c_ip_{n+i-1}(x)$.
Since the $y_i$'s are pairwise distinct, Identity~(\DC) implies that
$$q(x)=\prod _{i=1} ^{d}(x+x_i)
=\prod _{k=1} ^{e}(x+y_k)^{m_k}
$$ 
divides $g(x)$ as a polynomial in~$x$. Hence, there exists
another polynomial $h(x)$ such that
$$
g(x)=q(x)h(x).
$$

Now, by inspection, $g(x)$ is a polynomial of degree at most~$n+d-1$, while
$q(x)$ is a polynomial of degree~$d$. Hence, $h(x)$ is a polynomial of
degree at most~$n-1$, which therefore can be written as a linear
combination
$$
h(x)=\sum_{i=0}^{n-1}r_ip_i(x),
$$
for some constants $r_i\in\widehat K$, $i=0,1,\dots,n-1$.
By its definition, $g(x)$ is orthogonal to all $p_j(x)$ with $0\le
j\le n-1$. In other words, we have
$$
0=L\big(g(x)p_j(x)\big)
=\sum_{i=0}^{n-1}r_iL\big(q(x)p_i(x)p_j(x)\big),\quad 
\text{for }j=0,1,\dots,n-1.
$$
Equivalently, the rows of the matrix
$\big(L\big(q(x)p_i(x)p_j(x)\big)_{0\le i,j\le n-1}$ are linearly
dependent, and consequently the determinant on the right-hand side
of~(\DAa) vanishes, which in its turn implies that the determinant
on the left-hand side of~(\DAa), which is equal to the determinant
in (\DB), vanishes, as desired.

\medskip
Conversely, let the determinant in (\DB) be equal to zero, for some choice
of $x_i's$ in $\widehat K$. Again, without loss of generality
we assume that the first $m_1$ of the $x_i$'s
are equal to~$y_1$, the next $m_2$ of the~$x_i$'s are equal to~$y_2$,
\dots, and the last $m_e$ of the $x_i$'s are equal to~$y_e$, the
$y_j$'s being pairwise distinct. These assumptions imply again that
$$q(x)=\prod _{i=1} ^{d}(x+x_i)
=\prod _{k=1} ^{e}(x+y_k)^{m_k}.
$$ 
Using the equality
of Lemma~\TE, we see that the determinant on the right-hand side
of~(\DAa) must vanish. Thus, the rows of the matrix on this right-hand
side must be linearly dependent so that there exist
constants~$r_i\in\widehat K$, $i=0,1,\dots,n-1$, not all of them zero,
such that 
$$
0=\sum_{i=0}^{n-1}r_iL\big(q(x)p_i(x)p_j(x)\big),\quad 
\text{for }j=0,1,\dots,n-1.
$$
Now consider the polynomial $g(x)=\sum_{i=0}^{n-1}r_iq(x)p_i(x)$.
This is a non-zero polynomial of degree at most~$n+d-1$ which, by the last
identity, is orthogonal to $p_j(x)$, for $j=0,1,\dots,n-1$.
Hence there must exist constants~$c_i\in\widehat K$ such that
$$
g(x)=\sum_{i=1}^dc_ip_{n+i-1}(x).
$$
On the other hand, by the definition of $g(x)$, we have 
$$g^{(j-1)}(-y_k)=0,\quad 
\text{for }k=1,2,\dots,e\text{ and }j=1,2,\dots,m_k.
$$ 
We conclude that
$$
\sum_{i=1}^dc_ip^{(j-1)}_{n+i-1}(-x_k)=0,\quad 
\text{for }k=1,2,\dots,e\text{ and }j=1,2,\dots,m_k.
$$
This means that the rows of the matrix
$M_{m_1,m_2,\dots,m_e}(y_1,y_2,\dots,y_e)$ are linearly dependent, so
that, by Proposition~\TEa\ the polynomial
$R(x_1,x_2,\dots,x_d)$ in~(\DA) vanishes.\quad \quad \qed
\enddemo

We can now complete the proof of Theorem~\TA. 

\demo{Proof of Theorem \TA}
By Lemma~\TF\ we know that the symmetric polynomials\linebreak
$R(x_1,x_2,\dots,x_d)$ and $S(x_1,x_2,\dots,x_d)$ vanish only jointly.
If we are able to show that in addition both have the same highest
degree term then
they must be the same up to a multiplicative constant.
Indeed, the highest degree term in $R(x_1,x_2,\dots,x_d)$ is
obtained by selecting the highest degree term in each entry of the
matrix in the numerator of the fraction on the right-hand side
of~(\DA). Explicitly, this highest degree term is
$$
\frac {\displaystyle\det_{1\le i,j\le d}\left((-x_j)^{n+i-1}\right)}
{
\prod\limits _{1\le i<j\le d} ^{}(x_i-x_j)}
=
\prod _{j=1} ^{d}(-x_j)^n
\frac {\displaystyle\det_{1\le i,j\le d}\left((-x_j)^{i-1}\right)}
{
\prod\limits _{1\le i<j\le d} ^{}(x_i-x_j)}
=(-1)^{nd}
\prod _{j=1} ^{d}x_j^n,
$$
by the evaluation of the Vandermonde determinant.

On the other hand, the highest degree term in $S(x_1,x_2,\dots,x_d)$ is
obtained by selecting the highest degree term in each entry of the
matrix in the numerator of the fraction on the right-hand side
of~(\DB). Explicitly, this highest degree term is
$$
\det\limits_{0\le i,j\le n-1}\left(
\mu_{i+j}\prod _{\ell=1} ^{d}x_\ell\right)
=\prod _{j=1} ^{d}x_j^n
\det\limits_{0\le i,j\le n-1}\left(
\mu_{i+j}\right).
$$

Both observations combined, we see that
$$
C\frac {\displaystyle\det_{1\le i,j\le d}\left(p_{n+i-1}(-x_j)\right)}
{\prod\limits _{1\le i<j\le d} ^{}(x_i-x_j)}
=
\det\limits_{0\le i,j\le n-1}\left(
\mu^{i+j}\prod _{\ell=1} ^{d}(x_\ell+\mu)\right).
$$
with
$C=(-1)^{nd}\det_{0\le i,j\le n-1}(\mu_{i+j})$. This finishes the
proof of the theorem.\quad \quad \qed
\enddemo

\subhead 6. Fourth proof of Theorem \TA\ --- Heine's formula and Vandermonde 
determinants\endsubhead
The subject of this section is 
the random matrix-inspired proof of Theorem~\TA\
due to Br\'ezin and Hikami \cite{\BrHiAA}. Among all the proofs
of Theorem~\TA\ that I present in this paper, it is the only one that
does not need the knowledge of the formula beforehand. Rather, 
starting from the left-hand side, by the clever use of Heine's 
formula --- given in Lemma~\TFa\ below --- 
and some obvious manipulations, one arrives almost
effortlessly at the right-hand side.
The meaning of the formula in the context of random matrices has been
indicated in Subsection~2.3, and more specifically in~(\AZc) which
tells that the right-hand side of (\AA) can be interpreted as an
expectation of products of characteristic
polynomials of random Hermitian matrices. The random matrix flavour
of the calculations below is seen in the ubiquitous multivariate density
$$
\bigg(
\prod _{1\le i<j\le n} ^{}(u_i-u_j)\bigg)^2\,d\nu(u_1)d\nu(u_2)\cdots d\nu(u_n),
$$
which is, up to scaling, 
the density function for the eigenvalues of random Hermitian matrices.

\medskip
We prove (\AA) in the form
$$
{\det\limits_{0\le i,j\le n-1}\left(
\mu^{i+j}\prod _{\ell=1} ^{d}(x_\ell-\mu)\right)}
=
{\det\limits_{0\le i,j\le n-1}\left(\mu_{i+j}\right)}
\frac {\displaystyle\det_{1\le i,j\le d}\left(p_{n+i-1}(x_j)\right)}
{
\prod\limits _{1\le i<j\le d} ^{}(x_j-x_i)}.
\tag\YA
$$

As indicated, we need Heine's integral formula
(cf\. \cite{\SzegAA, Eq.~(2.2.11)} or
\cite{\IsmaAA, Cor.~2.1.3}).
Its proof is short enough that we provide it here for the sake of
completeness. The integral that appears in the formula can 
equally well be understood in the analytic or in  the formal sense.

\proclaim{Lemma \TFa}
Let $d\nu(u)$ be a density with moments $\nu_s=\int u^s\,d\nu(u)$,
$s=0,1,\dots$.
For all non-negative integers~$n$, we have
$$
\det_{0\le i,j\le n-1}(\nu_{i+j})=
\frac {1} {n!}\int \cdots\int\bigg(
\prod _{1\le i<j\le n} ^{}(u_i-u_j)\bigg)^2\,d\nu(u_1)d\nu(u_2)\cdots d\nu(u_n).
\tag\YB$$
\endproclaim

\demo{Proof}
We start with the left-hand side,
$$\align 
\det_{0\le i,j\le n-1}(\nu_{i+j})
&=\det_{0\le i,j\le n-1}\left(\int u_i^{i+j}\,d\nu(u_i)\right)\\
&=\int\cdots\int \bigg(
\prod _{i=0} ^{n-1}u_i^{i}\bigg)
\det_{0\le i,j\le n-1}\big(u_i^{j}\big)\,
d\nu(u_0)d\nu(u_1)\cdots d\nu(u_{n-1}).
\endalign$$
If in this expression we permute the $u_i$'s, then it remains
invariant, except for a sign that is created by the determinant. 
Let $\frak S_n$ denote the group
of permutations of $\{0,1,\dots,n-1\}$. If we average the above
multiple integral
over all possible permutations of the~$u_i$'s, then we obtain
$$\align 
\det_{0\le i,j\le n-1}(\nu_{i+j})
&=\frac {1} {n!}\int\cdots\int \sum_{\si\in\frak S_n}\sgn \si\,\bigg(
\prod _{i=0} ^{n-1}u_{\si(i)}^{i}\bigg)\\
&\kern3cm\cdot
\det_{0\le i,j\le n-1}\big(u_i^{j}\big)\,
d\nu(u_0)d\nu(u_1)\cdots d\nu(u_{n-1})\\
&=\frac {1} {n!}\int\cdots\int
\left(\det_{0\le i,j\le n-1}\big(u_i^{j}\big)\right)^2\,
d\nu(u_0)d\nu(u_1)\cdots d\nu(u_{n-1}).
\endalign$$
In view of the evaluation of the Vandermonde determinant, up to
a shift in the indices of the $u_i$'s, this is the right-hand side
of~(\YB).\quad \quad \qed
\enddemo

\demo{Proof of (\YA)}
Let $d\mu(u)$ be a density with moments $\mu_s=\int u^s\,d\mu(u)$,
$s=0,1,\dots$. We apply Lemma~\TFa\ with
$$
d\nu(u)=d\mu(u)\prod _{\ell=1} ^{d}(x_\ell-u).
$$
This yields
$$\multline
\det_{0\le i,j\le n-1}\left(
\mu^{i+j}\prod _{\ell=1} ^{d}(x_\ell-\mu)\right)\\
=
\frac {1} {n!}\int \cdots\int
\bigg(
\prod _{j=1} ^{n}\prod _{\ell=1} ^{d}(x_\ell-u_j)\bigg)
\bigg(
\prod _{1\le i<j\le n} ^{}(u_i-u_j)\bigg)^2\,d\mu(u_1)d\mu(u_2)\cdots d\mu(u_n).
\endmultline
\tag\YC$$
Let $\big(p_m(y)\big)_{m\ge0}$ be the sequence of monic orthogonal
polynomials with respect to the linear functional with
moments $\mu_s=\int u^s\,d\mu(u)$, $s=0,1,\dots.$
It is easy to see that
$$
\det_{1\le i,j\le N}\big(p_{j-1}(y_i)\big)
=
\det_{1\le i,j\le N}\big(y_i^{j-1}\big)
=\prod _{1\le i<j\le N} ^{}(y_j-y_i).
\tag\YD$$
We use this observation with $N=n+d$, where the role of the $y_i$'s
is taken by $u_1,u_2,\dots,u_n,\mathbreak x_1,x_2,\dots,x_d$.
Then (\YC) may be rewritten as
$$\multline
\det_{0\le i,j\le n-1}\left(
\mu^{i+j}\prod _{\ell=1} ^{d}(x_\ell-\mu)\right)\\
=
\frac {1} {n!\,\prod\limits_{1\le i<j\le d}(x_j-x_i)}\int \cdots\int
\det_{1\le i,j\le n+d}\pmatrix 
p_{j-1}(u_i),&1\le i\le n\hfill\\
p_{j-1}(x_{i-n}),&n+1\le i\le n+d
\endpmatrix\\
\cdot
\bigg(
\prod _{1\le i<j\le n} ^{}(u_j-u_i)\bigg)\,d\mu(u_1)d\mu(u_2)\cdots d\mu(u_n).
\endmultline
$$
By orthogonality of the $p_j(x)$'s, we have $\int u^sp_{j-1}(u)\,d\mu(u)=0$
for $0\le s<n\le j-1$. Since the expansion of the Vandermonde determinant
$\prod _{1\le i<j\le n} ^{}(u_j-u_i)$ consists entirely of monomials
$\pm u_1^{s_1}u_2^{s_2}\cdots u_n^{s_n}$ with $s_i<n$ for all~$i$, we
see that in the above expression we may replace the determinant in
the integrand by
$$\multline
\det_{1\le i,j\le n+d}\pmatrix 
p_0(u_1)&p_1(u_1)&\dots&p_{n-1}(u_1)&0&\dots&0\\
p_0(u_2)&p_1(u_2)&\dots&p_{n-1}(u_2)&0&\dots&0\\
\hdotsfor7\\
p_0(u_n)&p_1(u_n)&\dots&p_{n-1}(u_n)&0&\dots&0\\
p_0(x_1)&p_1(x_1)&\dots&p_{n-1}(x_1)&p_n(x_1)&\dots&p_{n+d-1}(x_1)\\
p_0(x_2)&p_1(x_2)&\dots&p_{n-1}(x_2)&p_n(x_2)&\dots&p_{n+d-1}(x_2)\\
\hdotsfor7\\
p_0(x_d)&p_1(x_d)&\dots&p_{n-1}(x_d)&p_n(x_d)&\dots&p_{n+d-1}(x_d)
\endpmatrix\\
=\det_{1\le i,j\le n}\big(p_{j-1}(u_i)\big)\cdot
\det_{1\le i,j\le d}\big(p_{n+j-1}(x_i)\big).
\endmultline
$$
If we substitute this above and use (\YD) again, we get
$$\multline
\det_{0\le i,j\le n-1}\left(
\mu^{i+j}\prod _{\ell=1} ^{d}(x_\ell-\mu)\right)\\
=
\frac {\det\limits_{1\le i,j\le d}\big(p_{n+j-1}(x_i)\big)} 
{n!\,\prod\limits_{1\le i<j\le d}(x_j-x_i)}\int \cdots\int
\bigg(
\prod _{1\le i<j\le n} ^{}(u_j-u_i)\bigg)^2\,d\mu(u_1)d\mu(u_2)\cdots d\mu(u_n).
\endmultline
$$
This gives indeed (\YA) once we apply Lemma~\TFa\ another time, now with 
$d\nu(u)=d\mu(u)$.\quad \quad \qed
\enddemo

\subhead 7. Hankel determinants of linear combinations of
Motzkin and Schr\"oder numbers\endsubhead
As described in Section~2, the origin of the author's discovery
of Theorem~\TA\ has been the interest in the evaluation of Hankel
determinants of linear combinations of combinatorial
sequences. Elouafi in \cite{\ElouAA} has the same motivation.
The point of (\AA) in this context is that
it provides a compact formula for $n\times n$ Hankel determinants
of a fixed linear combination of $d+1$ successive elements
of a (moment) sequence (the
left-hand side of~(\AA)) that does not ``grow" with~$n$. (The
right-hand side is a ``fixed-size" formula for fixed~$d$; the
dependence on~$n$ is in the index of the orthogonal polynomials.)
Elouafi provides numerous applications of Theorem~\TA\ to the evaluation
of Hankel determinants of linear combinations of Catalan, Motzkin, and
Schr\"oder numbers in \cite{\ElouAA, Sec.~3}. However, his treatment
of Motzkin and Schr\"oder numbers can be replaced by a better one. 

\medskip
The reader should recall that, if some of the $x_i$'s in~(\AA) should
be equal to each other, then on the right-hand side we would have to
use Proposition~\TEa\ in order to make sense of the right-hand side
of~(\AA). In order to apply the proposition, we must evaluate derivatives of
the orthogonal polynomials at~$-x_j$. How to accomplish this for
the orthogonal polynomials corresponding to Motzkin and Schr\"oder
numbers as moments by a recursive approach, 
is described in \cite{\ElouAA, Secs.~3.2 and~3.3}.
Here we show that one can be completely explicit.
We do not try to treat the most general case but rather restrict
ourselves to illustrate the main ideas by two examples.

\medskip
By Theorem~\TA\ with $x_i=0$ for all~$i$ and Proposition~\TEa\ with
$y_1=0$ and $m_1=d$, we have (see also \cite{\ElouAA, Eq.~(1.3)})
$$
\frac {\det\limits_{0\le i,j\le n-1}\left(
\mu_{i+j+d}\right)}
{\det\limits_{0\le i,j\le n-1}\left(\mu_{i+j}\right)}
=(-1)^{nd}
{\displaystyle\det_{1\le i,j\le
d}\left(\frac {p^{(j-1)}_{n+i-1}(0)} {(j-1)!}\right)}.
\tag\EA
$$

\medskip
We consider first the special case where
$\mu_n=M_n$, the {\it $n$-th Motzkin number}, defined by
$\sum_{n\ge0}M_n\,z^n=\frac {1-z-\sqrt{1-2z-3z^2}} {2z^2}$
(cf\. \cite{\StanBI, Ex.~6.37}). It is well-known (see \cite{\VienAE}
or \cite{\ElouAA, p.~1265}) that the associated orthogonal polynomials
satisfy the three-term recurrence~(\AB) with $s_i=t_i=1$ for all~$i$.
More explicitly, they
are $p_n(x)=U_n\left(\frac {x-1} {2}\right)$, where $U_n(x)$ is the
{\it $n$-th Chebyshev polynomial of the second kind}, which is defined by 
$$
U_n(x)=2xU_{n-1}(x)-U_{n-2}(x),
\tag\EAa$$
with generating function
$$
\sum_{n\ge0}U_n(x)z^n=\frac {1} {1-2xz+z^2}.
$$
From this generating function, we may now easily obtain an explicit
expression for $p^{(j)}_n(0)=\frac {d^j} {dx^j}U_n\left(\frac {x-1}
{2}\right)\big\vert_{x=0}$. Namely, we have 
$$
\sum_{n\ge0}\frac {d^j} {dx^j}U_n\left(\frac {x-1} {2}\right)z^n=
\frac {j!\,z^j} {\big(1-(x-1)z+z^2\big)^{j+1}}.
$$
Consequently, we get
$$
\align
\sum_{n\ge0}p^{(j)}_n(0)z^n&=
\frac {j!\,z^j} {\big(1+z+z^2\big)^{j+1}}\\
&=
\frac {j!\,z^j(1-z)^{j+1}} {\big(1-z^3\big)^{j+1}}\\
&=j!\,z^j\sum_{a=0}^{j+1}(-1)^a\binom {j+1}a z^a
\sum_{b\ge0}\binom {j+b}b z^{3b}.
\endalign
$$
Comparison of coefficients of $z^n$ then yields
$$
p^{(j)}_n(0)=
j!\,\sum_{b\ge0} 
(-1)^{n-j-b}\binom {j+1}{n-j-3b} 
\binom {j+b} {b}.
%\binom {j+b}b 
%n=j+a+3b 
%a=n-j-3b
$$
If we recall that $\det_{0\le i,j\le n-1}\left(M_{i+j}\right)=1$
(by~(\ABa) and the choice $t_i=1$ for all~$i$),
then we arrive at the identity
$$
{\det\limits_{0\le i,j\le n-1}\left(
M_{i+j+d}\right)}
=
 {\displaystyle\det_{1\le i,j\le
d}\big(A_{i,j}(n)\big)},
\tag\EB
$$
where
$$
A_{i,j}(n)=\sum_{b\ge0} 
(-1)^{b}\binom {j}{n+i-j-3b} 
\binom {j+b-1} {b}.
$$

\medskip
Second, we consider here the special case where
$\mu_n=r_n$, the {\it $n$-th (large) Schr\"oder number}, defined by
$\sum_{n\ge0}r_n\,z^n=\frac {1-z-\sqrt{1-6z+z^2}} {2z}$
(cf\. \cite{\StanBI, Second Problem on page~178}). 
Here, the associated orthogonal polynomials
satisfy the three-term recurrence~(\AB) with 
$s_0=2$, $s_i=3$ for $i\ge1$, and
$t_i=2$ for all~$i$ (see \cite{\CJKrAB, Case~(vii) in Sec.~4}).
More explicitly, their generating function satisfies
$$
\sum_{n\ge0}p_n(x)z^n=\frac {1+z} {1-(x-3)z+2z^2}.
$$
From this generating function, we may now obtain an explicit
expression for $p^{(j)}_n(0)$. Namely, we have 
$$
\sum_{n\ge0}\frac {d^j} {dx^j}p_n(x)z^n=
\frac {j!\,z^j(1+z)} {\big(1-(x-3)z+2z^2\big)^{j+1}}.
$$
Consequently, we get
$$
\align
\sum_{n\ge0}p^{(j)}_n(0)z^n&=
\frac {j!\,z^j(1+z)} {\big(1+3z+2z^2\big)^{j+1}}\\
&=\frac {j!\,z^j} {(1+z)^{j}(1+2z)^{j+1}}\\
&=j!\,z^j\sum_{a=1}^j\frac {c_a} {(1+z)^a}
+j!\,z^j\sum_{b=1}^{j+1}\frac {d_b} {(1+2z)^b},
\endalign
$$
where
$$
c_a=\frac {1} {(j-a)!}
\frac {d^{j-a}} {dz^{j-a}}\frac {1} {(1+2z)^{j+1}}\bigg\vert_{z=-1}
=
(-1)^{j+1}2^{j-a}\binom {2j-a}j
%\frac {(-1)^{j-a+a}2^{j-a}(j+1)(j+2)\cdots(2j-a)} {(j-a)!}
$$
and
$$
d_b=\frac {1} {(j+1-b)!}
\frac {d^{j+1-b}} {dz^{j+1-b}}\frac {1} {(1+z)^{j}}\bigg\vert_{z=-1/2}
=
(-1)^{j+b+1}2^{j}\binom {2j-b}{j-1}.
%\frac {(-1)^{j+b+1}2^{2j+1-b}(j)(j+1)\cdots(2j-b)} {(j+1-b)!}
$$
Then, by comparing coefficients of $z^n$, we obtain
$$
\multline
p^{(j)}_n(0)
=j!\sum_{a=1}^j(-1)^{n+1}2^{j-a}\binom {2j-a}j 
\binom {a+n-j-1}{n-j}\\
+j!\sum_{b=1}^{j+1}(-1)^{n+b+1}2^{n}\binom {2j-b}{j-1} 
\binom {b+n-j-1}{n-j}.
\endmultline
$$
If we recall that $\det_{0\le i,j\le
n-1}\left(r_{i+j}\right)=2^{\binom n2}$ 
(by~(\ABa) and the choice $t_i=2$ for all~$i$),
%or cf\. e.g\. \cite{\EuFuAA, Cor.~2.5(i)}), 
then we arrive at the identity
$$
{\det\limits_{0\le i,j\le n-1}\left(
r_{i+j+d}\right)}
=(-1)^{\binom {d+1}2}2^{\binom n2}
 {\displaystyle\det_{1\le i,j\le
d}\big(B_{i,j}(n)\big)},
\tag\EC
$$
where
$$
\multline
B_{i,j}(n)=
\sum_{a=1}^{j-1}2^{j-a-1}\binom {2j-a-2}{j-1} 
\binom {a+n+i-j-1}{n+i-j}\\
+\sum_{b=1}^{j}(-1)^{b}2^{n+i-1}\binom {2j-b-2}{j-2} 
\binom {b+n+i-j-1}{n+i-j}.
\endmultline$$

\medskip
As a final remark, we point out that the above treatment of the
``Motzkin case" is one which only applies in a specific situation,
whereas the above treatment of the ``Schr\"oder case" works
for all moment sequences for orthogonal polynomials which are
generated by a three-term recurrence~(\AB) where the coefficient
sequences $(s_i)_{i\ge0}$ and $(t_i)_{i\ge0}$ become constant
eventually, that is,
where $s_i\equiv s$ and
$t_i\equiv t$ for large enough~$i$.
For, under this assumption, the generating function
$\sum_{n\ge0}p_n(x)z^n$ for the orthogonal polynomials
$\big(p_n(x)\big)_{n\ge0}$ is a rational function with denominator
$1-(x-s)z+tz^2$, a quadratic polynomial, as in the special case of
Schr\"oder numbers, where we had $s=3$ and $t=2$.

\subhead 8. Recursiveness of Hankel determinants of linear
combinations of moments\endsubhead
In \cite{\DoFSAA, Conj.~5},
Dougherty, French, Saderholm and Qian conjectured that
the Hankel determinants of a linear combination of {\it Catalan
numbers} $C_n=\frac {1} {n+1}\binom {2n}n$,
$$
\det_{0\le i,j\le n-1}
\left(\sum_{k=0}^d\la_kC_{i+j+k}\right),
$$
satisfy a linear recurrence with constant coefficients of order
$2^d$. This conjecture was proved by Elouafi in \cite{\ElouAA,
Theorem~6} by using Theorem~\TA\ and special properties of the
orthogonal polynomials corresponding to the Catalan numbers as moments. 
However, Theorem~\TA\ implies a much more general result. 
This is what we make explicit in the corollary below.
For the statement of the corollary, we need to recall the definition
of the elementary symmetric polynomials
$$
e_k(x_1,x_2,\dots,x_d):=\sum_{1\le i_1<i_2<\dots<i_k\le d}x_{i_1}x_{i_2}
\cdots x_{i_k}.
$$

\proclaim{Corollary \TG}
Within the setup in Section~{\rm1}, let $s_i\equiv s$ and $t_i\equiv t$ for
$i\ge1$. Furthermore, let
$$
H_n=t_0^{-(n-1)}t^{-\binom {n-1}2}
\det_{0\le i,j\le n-1}\Bigg(\sum _{k=0} ^{d}\la_k\mu_{i+j+k}\Bigg),
\tag\ZA$$
where the $\la_k$'s are some constants\footnote{We may also think of
the $\la_k$'s as variables.} in the ground field~$K$, 
$k=0,1,\dots,d-1$, and $\la_d=1$.
Then the sequence $(H_n)_{n\ge0}$ of (scaled) Hankel determinants of linear
combinations of moments
satisfies a linear recurrence of the form
$$\sum_{i=0}^{2^d}c_iH_{n-i}=0,\quad \text{for }n> 2^d,
\tag\ZB$$
for some constants $c_i\in K$, normalised by $c_0=1$. 
Explicitly, these constants can be computed as the
coefficients of the
characteristic polynomial (in~$x$) $\sum_{i=0}^{2^d}c_ix^{2^d-i}$
of the tensor product of\/ $2\times2$ matrices
$$
\pmatrix x_1+s&t\\-1&0\endpmatrix
\otimes
\pmatrix x_2+s&t\\-1&0\endpmatrix
\otimes
\dots
\otimes
\pmatrix x_d+s&t\\-1&0\endpmatrix,
\tag\ZBa$$
where $\la_k=e_{d-k}(x_1,x_2,\dots,x_d)$.
In particular, we have
$$
c_1=-\sum_{j=0}^d\la_js^j
\tag\ZC$$
and
$$
c_{2^d-i}=t^{d(2^{d-1}-i)}c_i,\quad \text{for }i=0,1,\dots,2^d.
\tag\ZD$$
\endproclaim

\remark{Remarks}
(1) The reader should note that, by (\ABa), the scaling in (\ZA) is exactly the
value of the determinant
$$
\det_{0\le i,j\le n-1}\big(\mu_{i+j}\big)
$$ 
in the denominator on the left-hand side of (\AA). 

\medskip
(2) A small detail is that the proof of Corollary~\TG\ given below shows
that, if $t_i\equiv t$ for {\it all\/}~$i$, then the recurrence~(\ZB)
holds even for $n\ge 2^d$. 

Furthermore, an inspection of
the proof of the corollary shows that, 
if $s_i\equiv s$ for $i>N$ and $t_i\equiv t$ for $i\ge N$, where~$N$ is some
positive integer, then 
the recurrence~(\ZB) still holds, but only for $n\ge 2^d+N$.

\medskip
(3) The formula (\ZC) for the coefficient $c_1$ 
is a far-reaching generalisation 
of Conjecture~6 in \cite{\DoFSAA}, while the symmetry relation~(\ZD)
is a far-reaching generalisation 
of Conjecture~7 in \cite{\DoFSAA}.

\medskip
(4) In view of the specialisations listed in items (i)--(xi)
and (xiv)--(xviii) in the list given in \cite{\CJKrAB, Sec.~4},
Corollary~\TG\ implies that the (properly scaled) Hankel determinants
of linear combinations of numerous combinatorial sequences satisfy
a linear recurrence with constant coefficients, which, aside from
the already mentioned Catalan numbers, include Motzkin numbers,
central binomial coefficients, central trinomial coefficients,
central Delannoy numbers, Schr\"oder numbers, Riordan numbers,
and Fine numbers.
\endremark

\demo{Proof of Corollary \TG}
We use Theorem~\TA\ in the equivalent form~(\BB).
In order to apply this identity, we write
$$
\sum_{j=0}^d\la_jx^j=
\prod _{i=1} ^{d}(x+x_i),
$$
with the $x_i$'s in the algebraic closure of our ground
field~$K$. Equivalently, 
$$\la_k=e_{d-k}(x_1,x_2,\dots,x_d).$$
For the moment, we assume that the $x_i$'s are pairwise distinct in
order to avoid a zero denominator in~(\BB).
We will get rid of this restriction in the end by a limiting argument.
(An alternative would be to base the arguments on Proposition~\TEa.
This would however be more complicated.)

What
(\BB) affords is to express $H_n$ in terms of a $d\times d$
determinant with entries $f_{n+i-1}(x_j)$, $1\le i,j\le d$.
If we expand the determinant on the right-hand side of~(\BB) according
to its definition, then we obtain a
linear combination of products of the form
$
\prod _{i=1} ^{d}f_{n+\ta_i}(x_i)$
with constant coefficients,
where $\ta_i\in\{0,1,\dots,d-1\}$ for $i=1,2,\dots,d$. 
Now, by (\BA), for each fixed~$i$ and~$\ta_i$ 
the sequence $\big(f_{n+\ta_i}(x_i)\big)_{n\ge0}$
satisfies the recurrence relation
$$
g_{n}-(x_i+s)g_{n-1}+tg_{n-2}=0,\quad \text{for }n\ge3.
\tag\ZE
$$
Hence, each product sequence 
$\left(\prod _{i=1} ^{d}f_{n+\ta_i}(x_i)\right)_{n\ge0}$
satisfies the same recurrence relation,
namely the one resulting from the (Hadamard) product
of the recurrences (\ZE) over $i=1,2,\dots,d$.
From the proof of \cite{\StanBI, Theorem~6.4.9}
(which is actually a much more general theorem), it follows
immediately that the order of
this ``product'' recurrence is at most $2^d$.

\medskip
In order to obtain the explicit description of the ``product"
recurrence in the statement of the corollary, we have to recall
the basics of the solution theory of (homogeneous) linear recurrences
with constant coefficients. This theory says that one has to
determine the zeroes of the characteristic polynomial of the
recurrence; the powers $\al^n$
of the zeroes~$\al$ multiplied by powers $n^e$ of~$n$, where the
exponent~$e$ is less than the multiplicity of~$\al$,
form a basis of the solution space of the recurrence.

The characteristic polynomial of the recurrence~(\ZE) is
$$
x^2-(x_i+s)x+t,
\tag\ZF$$
which is also equal to the characteristic polynomial of the $2\times
2$ matrix
$$
\pmatrix x_i+s&t\\-1&0\endpmatrix.
\tag\ZG$$
Let $y_{i,1}$ and $y_{i,2}$ be the zeroes of the
polynomial~(\ZF). Then the powers $\big(y_{i,1}^n\big)_{n\ge0}$ and
$\big(y_{i,2}^n\big)_{n\ge0}$ form a basis of solutions to the
recurrence~(\ZE) if $y_{i,1}$ and $y_{i,2}$ are distinct, while
otherwise the sequences $\big(y_{i,1}^n\big)_{n\ge0}$ and
$\big(ny_{i,1}^n\big)_{n\ge0}$ form a basis. 
The product recurrence that we want to find
is one for which all products $\left(
n^{e-d}\prod _{i=1} ^{d}y_{i,\ep_i}^n\right)_{n\ge0}$ are solutions,
for $\ep_i\in\{1,2\}$ and~$e$ bounded above by the sum of the
multiplicities of the~$y_{i,1}$'s. 
Equivalently, the characteristic polynomial of
the desired product recurrence is one which has all products
$\prod _{i=1} ^{d}y_{i,\ep_i}$,
where $\ep_i\in\{1,2\}$, as zeroes, with multiplicities equal
to the sum of the multiplicities of the $y_{i,1}$'s minus~$d-1$.
It is a simple fact of linear algebra that such a polynomial 
is the characteristic polynomial of the tensor
product of all matrices~(\ZG), that is, the matrix in~(\ZBa).
(For, the eigenvalues of the tensor product of matrices are
all products of eigenvalues of the individual matrices.)
This proves the assertion about the explicit form of the
recurrence~(\ZB) in the case of pairwise distinct~$x_i$'s.

Since everything --- the expressions in~(\BB), the coefficients of
the recurrence~(\ZB) --- is polynomial in the $x_i$'s, we may finally drop
that restriction.

\medskip
The coefficient $c_1$ is the coefficient of $x^{2^d-1}$ in the
characteristic polynomial of (\ZBa). It is easy to see that this is
$$
-
\prod _{i=1} ^{d}(x_i+s)=-\sum_{j=0}^de_{d-j}(x_1,x_2,\dots,x_d)s^j
=-\sum_{j=0}^d\la_js^j,
$$
proving (\ZC).

\medskip
Finally, the symmetry relation (\ZD) is a consequence of the inherent
symmetry of a linear recurrence of order~2. In order to make this
visible, let $\widehat f_n(x)=t^{-n/2}f_n(x)$. Then, from~(\BA) we see
that 
$$
\widehat f_n(x)-(x+s)t^{-1/2}\widehat f_{n-1}(x)+\widehat f_{n-2}(x)=0,
\quad \text{for }n\ge3.
\tag\ZH$$
Now, a recurrence can be read in the forward direction --- that is,
we compute the $n$-th term of the sequence from lower order terms ---
or in the backward direction --- that is, we compute the $n$-th term
from higher order terms. In this sense we see that the
recurrence~(\ZH) is the same regardless whether it is read in forward
or backward direction. Consequently, the recurrence relation for
the Hadamard product $
\prod _{i=1} ^{d}\widehat f_n(x_i)$ must also be symmetric, that is,
the same regardless whether it is read in forward or in backward
direction. If one then substitutes back $\widehat
f_n(x)=t^{-n/2}f_n(x)$ in that symmetric recurrence, the
relation~(\ZD) follows.\quad \quad \qed
\enddemo

\subhead Acknowledgements\endsubhead
I thank Mourad Ismail for insightful discussions and Arno Kuijlaars
for bringing \cite{\BaDSAA} to my attention.

\enddocument